\documentclass[11pt, 3p]{article}
\usepackage[a4paper, total={6in, 8in}]{geometry}
\usepackage[utf8]{inputenc}
\usepackage[english]{babel}
\usepackage{amsmath}
\usepackage{amsthm}
\usepackage{amsfonts}
\usepackage{amssymb}
\usepackage{bm}
\usepackage{mathtools}
\usepackage{csvsimple}
\usepackage[backend=biber,style=numeric,autocite=plain,sorting=nyt,maxbibnames=99]{biblatex}
\usepackage{csquotes}
\usepackage{color}
\usepackage{algorithm}
\usepackage{algpseudocode}
\usepackage{xcolor}
\usepackage{tikz}
\usepackage{algorithmicx}
\usepackage{algpseudocode}
\usepackage{comment}
\usepackage{setspace}
\usepackage{graphicx}
\usepackage{subcaption}

\newtheorem{definition}{Definition}

\newcommand\arrowHeadScale{1.5}

\title{Optimization for Evaluating the Practical Capacity of a Transshipment Yard}
\date{}

\begin{document}

\maketitle

\begin{center}

{\large
Anna Russo Russo $^{a}$, Roberto Mancini $^{b}$, Gianpaolo Oriolo $^{a}$, \\ Veronica Piccialli $^{c}$, Davide Ussai $^{b}$
}

\vspace{3ex}

\vspace{2ex}$^{a}$ Dipartimento di Ingegneria Civile e Ingegneria Informatica, \\Università di Roma Tor Vergata, Italy \\\texttt{anna.russo.russo@uniroma2.it, oriolo@disp.uniroma2.it}
 \\ 
\vspace{2ex}$^{b}$ Rete Ferroviaria Italiana RFI, Italy \\\texttt{r.mancini@rfi.it, d.ussai@rfi.it}
 \\ 
\vspace{2ex}$^{c}$ Dipartimento di ingegneria Informatica, Automatica e Gestionale, Università di Roma La Sapienza, Italy \\\texttt{veronica.piccialli@uniroma1.it}
\end{center}

\newpage

\begin{abstract}
In order to increase rail freight transportation in Italy, Rete Ferroviaria Italiana (RFI) the Italian railway infrastructure manager, is carrying out several investment plans to enhance the Transshipment Yards, that act as an interface between the rail and road networks. The need is to increase their practical capacity, i.e. the maximum number of train services that can be inserted without altering the current timetable while respecting all relevant constraints. Several factors influence the practical capacity of a transshipment yard: physical resources (such as tracks and vehicles for loading/unloading); constraints on the possible time slots of individual operations; constraints on the length of time a train must stay in the yard, that follow from both timetable requirements that are settled by the (prevalent) main line and from administrative and organisational issues in the yard. In this paper, we propose a {\sc milp}-based optimization model that is based on the solution of a suitable {\em saturation} problem, that deals with all these constraints and that can be used for evaluating the practical capacity of a transshipment yard both in its current configuration and in any plausible future configuration. The model provides operational details, such as routes and schedules, for each train service, and allows to impose {\em periodic} timetables and schedules that keep the daily management of the yard easier. Both the model and its solutions are validated on a real Italian transshipment yard, located at Marzaglia, on different scenarios corresponding to different investment plans of RFI. The results show that proper investments allow to get a feasible timetable with a period of 24 hours with doubles the number of current train services.
\end{abstract}

\section{Introduction}

\bigskip

Global rail freight transportation has been interested in a steady growth over the past two decades \cite{report_rail_2019} and the same trend has been observed in Italy.
This growth in rail freight transportation demand has posed the necessity to  strengthen and extend the so-called \textit{Last Mile} area of the network, as many Italian transshipment yards are becoming more and more congested.
Transshipment yards serve as an interface between roads and rail networks in the intermodal transportation of freight.
Freight trains access these yards through an arrival and departure station that is connected to the main line, then are moved inside the internal infrastructure of the yard, where load units are moved from trains onto road vehicles and vice versa. Later, the trains move again towards the arrival and departure station and finally leave the yard at a scheduled departure time.
Both the timing and the physical details of these movements are highly dependent on the resource competition with the other freight trains.
In order to be able to allocate slots for new trains services, Rete Ferroviaria Italiana (RFI), the Italian railway infrastructure manager, is considering and carrying out several investment plans for its transshipment yards.
However, rail infrastructure is extremely expensive and has a huge impact on the surrounding environment.
In this context, having a tool for the quantitative evaluation of a terminal capacity, both in the actual or in any hypothetical configuration, is of high interest for RFI, as that would facilitate a data-driven comparison between different investment scenarios. In the literature, the capacity of a railway network has been widely studied, and multiple definitions can be found. \cite{abril_2008} studies the capacity of a single-line railway network, suggesting a distinction between the \textit{theoretical} and the \textit{practical} capacity. This distinction can be extended to transshipment yards, where the theoretical capacity can be intended as the maximum number of trains that can be routed safely in a transshipment yard under ideal conditions, like uniform arrival times, and suitable departure times.
On the other hand, the practical capacity is the maximum number of trains that can be served by the transshipment yard, under realistic conditions and at a reasonable level of reliability. 

The theoretical capacity is virtually unachievable, therefore our interest is in estimating the practical capacity of a transshipment yards, starting from realistic assumption, e.g. dealing with timetable that is not ideal but is periodic, that is what happens in practice.

In a real timetable for an Italian transshipment yard arrivals and departures are usually concentrated in just a fraction of the day so that freight trains do not conflict with maintenance hours and passengers trains on the line. This implies that trains usually have to wait a long time inside the transshipment yard after they have been processed, occupying precious resources.

Summarizing, several factors contribute to influencing the practical capacity of a freight yard: physical resources (such as tracks and vehicles for loading/unloading), constraints on the possible time slots of individual operations, constraints on the length of time a train must stay in the yard, that follow from both timetable requirements that are settled by the (prevalent) main line and from administrative and organisational issues in the yard (due to train crews turnover, for example).

In this paper we propose an optimization-based approach able to deal with all these constraints. We use our approach for a twofold purpose: evaluating the practical capacity of a transshipment yard both in the current configuration and in any plausible future configuration, and provide operational instructions to use the yard at full capacity.

Furthermore, our model allows to impose {\em periodic} timetables and schedules, i.e., leading to the design of periodic solutions which are quite desirable in order to keep the daily management of the yard easier.

The rest of the paper is organized as follows: in Section~\ref{section:rfi_literature} we provide a brief overview of the existing methods for the capacity evaluation in the context of rail. In Section~\ref{section:problem_description} we give a formal definition of the optimization problems that we propose for the capacity analysis of a transshipment yard, that we name as the \textit{Yard Saturation Problem} ({\sc ysp}) and the \textit{Periodic Yard Saturation Problem} ({\sc pysp}). In Section~\ref{section:disjunctive_graphs} we model both these problem through suitably disjunctive graphs. In Section~\ref{section:rfi_milp} we define a Mixed Integer Linear Programming formulation for {\sc pysp} based on those graphs. In Section~\ref{section:rfi_results} we validate our approach by some computational tests carried out on a real transshipment yard, located at Marzaglia, Italy.

\section{Literature Review}
\label{section:rfi_literature}

Literature on capacity evaluation in the context of rail transportation comprises contributions for different systems, such as railway lines, stations, junctions and transshipment yards. More in general, in spite of many hybrid contributions, we may classify these works according to the methodology that is used to solve the model. 
We distinguish between \textit{simulation-based} methods, \textit{analytical} methods and \textit{optimization-based} methods.

Simulation-based methods \cite{simulation_reis_2019, coviello_2015, coviello_2019} reproduce the behavior of a railway system to estimate the practical capacity, usually relying on commercial software.
Analytical methods \cite{analytical_kraft_1988, analytical_burdett_2006, burdett_2015, analytical_bevrani_2015, mussone_2013} rely on closed mathematical formulas and generally provide a good estimate of the theoretical capacity of a system.
Notable examples of analytical methods are those  elaborated by the International Union of Railways (UIC) and proposed in \cite{uic_1983} and \cite{uic_2013}.
Optimization-based methods consist in estimating the practical capacity by solving an optimization problem. Most of the times, this optimization problem closely resembles a scheduling problem, since the evaluation of the maximum number of trains that can use a railway system usually involves finding a feasible schedule for the movements of such set of trains.

The most popular problem employed in optimization-based methods is the \textit{saturation problem}. At a very high level, the problem consists of determining the maximum number of trains that it is possible to serve in a given system, such as a railway line, a station, a junction, or a transshipment yard. The problem was likely mentioned for the first time in \cite{delorme_2001}, but since then, the problem has been declined in numerous variants.
\cite{delorme_2001} considers the problem of scheduling the maximum number of trains inside a junction, starting from an empty timetable. For each train, the model has to choose the most convenient route and departure time, without introducing waiting times inside the junction. Two different heuristics are proposed for its solution.
In \cite{ingolotti_2004} the authors focus on the problem of adding train services, from a given set of candidates to a current fixed timetable in the context of a railway line. The problem is again solved through a heuristic algorithm.
The goal of the saturation problem addressed in \cite{burdett_2009} is inserting additional train services to a current timetable with no or minimal disruption of the current schedule for a railway line. This problem is formulated as a job-shop scheduling with time-window constraints, described by means of a disjunctive graph, and then solved heuristically through a sequence of scheduling problems where many decision variables are fixed.
\cite{cacchiani_2010} deals with the problem of adding extra freight trains to a current timetable for a railway network, minimizing the shift with respect to a given suggested schedule for the additional trains. The problem is formulated through a MIP model and is solved with a Lagrangian heuristic. The previous approach is refined in \cite{cacchiani_2016} where some operational constraints, such as the maximum number of trains that a junction or station can host at a given time, are taken into account.
\cite{jiang_2014} considers the problem of scheduling additional train services on a transit line, with respect to the current timetable and schedule, while minimizing both the travel times of the new train services and the disruption of the schedule. The problem is then solved through mixed integer linear programming.
\cite{pellegrini_2017} defines a multi-objective {\sc milp} model for maximizing both the numbers of additional passenger and freight trains scheduled on a railway line; the model is then solved by the $\epsilon$-constraint method. In this approach, the pre-existing timetable is considered to be fixed and the only decision left to the optimization model is which saturating trains, out of a large candidates set, can be inserted into this timetable without incurring in conflicts. The problem is solved through a heuristic decomposition of the complete time-horizon in smaller windows.
The methodology proposed in \cite{pellegrini_2017} was taken as a starting point for the work presented in \cite{pascariu_2021}, where the same ILP model is adapted to the analysis of an Italian freight node, composed of the passenger station of Novara Centrale and the transshipment yard of Novara Boschetto.
The model is then solved with the same saturation approach inside a novel saturation strategy implementing non-uniform priority patterns in the saturating trains set.
The issue of saturating a timetable with heterogeneous candidate train services is addressed also in \cite{liao_2024}, for saturating a macroscopic model of a railway network. Here the problem is modeled through a time-space network and the Pareto front of saturated timetables is obtained with respect to different train services categories.

\subsection{Our Contribution}
We propose an optimization model to evaluate the practical capacity of a transshipment yard. The idea is, given a certain layout and a fixed set of train services with their timetable, to evaluate how many slots can be added to the timetable to serve new candidate train services in reasonable arrival and departure time windows. Both the model and its solutions are validated on a real transshipment yard, located at Marzaglia, Italy.
The model builds upon a suitable {\em saturation problem}, where we want to insert additional train services on top of the current services with no disruption of their timetable. Our saturation problem is then solved through (exact) Integer Programming techniques. As discussed in the previous section, several models have been proposed so far for the solution of saturation problems, we therefore discuss in the following which, to the best of our knowledge, are the main novelties of our approach. 

First, our model is microscopic, as it provides operational details, such as routes and schedules, for each train service, being able to deal with the complexity of transshipment yards, where a large degree of parallelism allows to define, for a same train service, different routes and operations plans.  Moreover, our model is flexible enough to allow to evaluate different scenarios involving either infrastructural interventions in the transshipment layout, or changes in the shift duration and/or operational choices (as an example consider different vehicles for loading/unloading trains). This flexibility was particularly important for RFI since it allowed to evaluate on Marzaglia the effectiveness in terms of capacity improvement of disjunctive investment plans with different costs.

Finally, our model allows to impose {\em periodic} timetables and schedules, i.e., leading to the design of periodic solutions. While that is rather challenging from both a modelling and a computational point of view, periodic solutions are quite desirable in order to keep the daily management of the yard easier. We point out that, when solving a saturation problem, it is indeed sensible to evaluate whether, simply adding more replicas of current train services with lower frequency, it is possible to increase the regularity of the timetable. Once again this has been validated at Marzaglia where it was shown that, without any structural intervention, it is possible to add new train services as to increase the number of weekly trains from 47 to 72 and reduce the period from a week to 24 hours.

\section{The Yard Saturation Problem}
\label{section:problem_description}

In order to evaluate the practical capacity of a transshipment yard we deal with a suitable saturation problem. The basic problem will be introduced in Section~\ref{section:model_intro} and formalized in Section~\ref{section:formaldef}. We will then extend it to a periodic setting in Section~\ref{section:periodic_scheduling}. 

Our aim is twofold: on the one hand, we aim at presenting saturation problems that are interesting on their own, and that we hope might be relevant for other railway systems; on the other, we want to validate our models and solutions on a real yard that is located at Marzaglia, Italy. We therefore find convenient to start with a discussion about the main features of a transshipment yard.

\subsection{Transshipment yards}\label{sc:trans}
\begin{figure}
\centering
\includegraphics[scale=0.4]{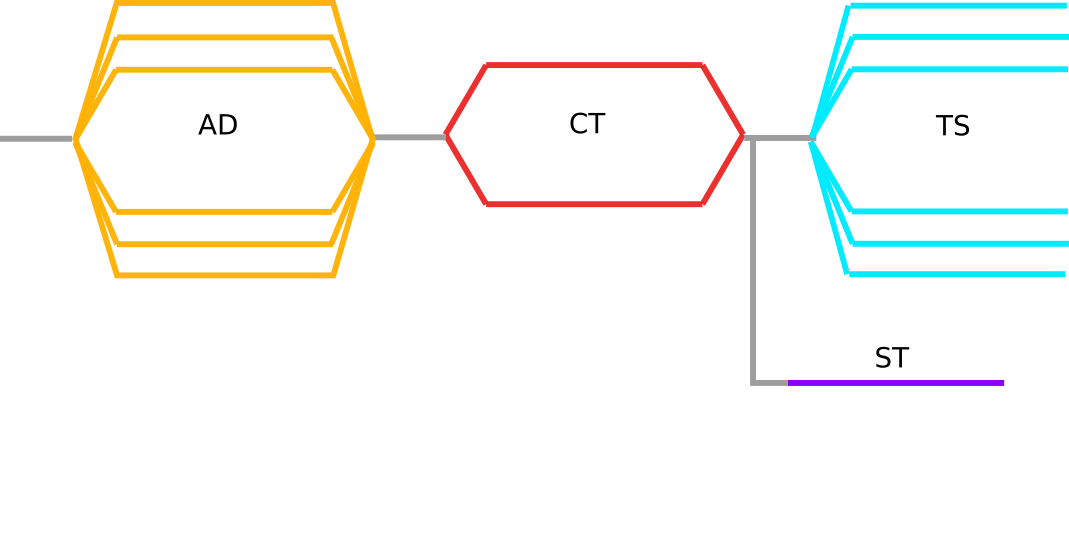} 
\caption{The layout of the transshipment yard of Marzaglia. In the arrival/departure area ($AD$) there are 6 tracks as well as in the transshipment area ($TS$); moreover the side tracks area ($ST$) is equipped with one track. Finally there are 2 tracks ($CT$) connecting the arrival/departure area to the transshipment area.}
\label{fig:marzaglia_layout}
\end{figure}

The layout of a transshipment yard is usually quite simple. There are different areas, each equipped with one or more tracks. For our purposes, relevant areas are: the \textit{arrival/departure area}, that serves as a buffer between the yard and the main rail network and usually hosts many parallel tracks; the \textit{transshipment area}, an area where trains are placed during the unloading and loading operations and that is equipped with parallel tracks and {\em vehicles} such as reach stackers, forklifts and rail-mounted cranes; possibly, the \textit{side tracks area}, where no operation takes place but trains are allowed to wait. Finally, \textit{connecting tracks} may connect the different areas. We provide a sketch for the layout of the yard in Marzaglia in Figure~\ref{fig:marzaglia_layout}, with an arrival/departure area with 6 tracks, a transshipment area with 6 tracks, a side tracks area with one track and 2 connection tracks.
 
In a transshipment yard each freight train undergoes a sequence of operations that usually involve: the arrival on a track of the arrival/departure area; the unloading and loading of the goods, performed by one or more vehicles in transshipment area; the return to the arrival/departure area; the  departure from the arrival/departure area.

Note that, while most of the above operations are mandatory for each train in the yard, there are usually {\em different plans} that are available for a same train: e.g. two plans for a train in the Marzaglia Yard may differ because one plan schedules a stop in the side tracks area and the other plan does not. Moreover, a same plan may corresponds to different routes because a same operation can be usually taken on parallel tracks (see Figure~\ref{fig:marzaglia_layout}). The formal definition of a {\em plan} will be given in Section~\ref{section:formaldef}.

\subsection{A model to evaluate the practical capacity of a transshipment yard}
\label{section:model_intro}
In order to evaluate the practical capacity of a transshipment yard we propose the solution of some saturation problems. In this section we start with the model, the {\em Yard Saturation Problem} ({\sc yps}). 

In the {\sc yps} we assume that we are given a set $T$ of {\em train services} that are operated at some (transshipment) yard together with their current {\em timetable}, i.e., the set of  arrival and departure times. We provide in Table~\ref{table:current_timetable} the current timetable for arrivals in Marzaglia: it is a weekly timetable comprising of {\em 47 trains} with {\em 11 train services} that are operated either daily, or three times a week, or twice a week. 

\begin{table}
\centering
\begin{tabular}{|c|c|c|c|c|c|c|}
\hline
Train & Mon & Tue & Wed & Thu & Fri & Sat \\
\hline
$t_1$ &    00:14 & 00:14 & 00:14 & 00:14 & 00:14 & 00:14 \\
$t_2$ &    & 04:52 & & & & 04:52 \\
$t_3$ &    & & 06:46 & & 06:46 & \\
$t_4$ &   09:11 &09:11 &09:11 &09:11 &09:11 &09:11 \\
$t_5$ &    &14:47 & &14:47 & & \\
$t_6$ &   15:54 &15:54 &15:54 &15:54 &15:54 &15:54 \\
$t_7$ &    & &16:10 & &16:10 & \\
$t_8$ &    &17:38 & &17:38 & &17:38 \\
$t_9$ &   20:38 &20:38 &20:38 &20:38 &20:38 &20:38 \\
$t_{10}$ &20:40 &20:40 &20:40 &20:40 &20:40 &20:40 \\
$t_{11}$ &20:42 &20:42 &20:42 &20:42 &20:42 &20:42 \\
\hline
\end{tabular}
\caption{Current timetable at Marzaglia for arrival. A similar timetable exists for departure.}
\label{table:current_timetable}
\end{table}

We will later deal with questions related to periodic timetables and saturation problems, see Section~\ref{section:periodic_scheduling}, where we will consider that each train service in $T$ is operated regularly with a certain given frequency. However, for the the {\sc yps} we assume that each train service in $T$ is operated just once.

Besides $T$ and the timetable for each train service in $T$, we are also given a set of candidate additional train service $T'$. For each train service $t\in T'$ we are also given time-windows, usually of size 60 mins, for both the arrival and the departure time (we motivate this below). We then would like to find the largest subset $T^*\subseteq T'$ of train services that can be served on top of the current set $T$. However, in doing so, there are a few constraints to meet. 

Trivially, we should meet any constraints related to plans, capacities etc. Assume that these constraints are satisfied -- we will provide formal definitions in Section~\ref{section:formaldef} -- and therefore we have a schedule and a timetable for each train in $T^*\cup T$. We require that the new timetable is such that there are no disruptions for current services, i.e., the train services in $T$ must keep the same arrival and departure times as in the current timetable. As for the train services in $T^*$, if any, arrival and departure times must be chosen within the given time-windows. 

We point out that, while arrival and departure times are fixed as described above, we are free to choose any plan, schedule and route for any $t\in T\cup T^*$: in particular, we are free to choose for train services in $T$ plans and schedules different from the current ones, provided that we do not alter the current arrival and departure times. 

Before moving to a more formal description of the {\sc yps}, we go back to the issue of time windows for arrival and departure times for trains in $T'$. Upper bounds are not uncommon in the literature, but lower bounds are. In our context, lower bounds are justified by the fact that, in practice, trains can only move from the transshipment yard according to a timetable which regulates movements on the main line. Acknowledging these constraints is {\em essential} for the capacity evaluation of a transshipment yard, also because the length of the stay of a train in the transshipment yard has a huge impact on the capacity of the latter.

\subsection{A formal definition of the Yard Saturation Problem}\label{section:formaldef}

In the {\em Yard Saturation Problem} we are given:
\begin{itemize}
	\item a set $S$ of resources (e.g. arrival/departure tracks), each with a capacity $u_s \in \mathbb{N}_+$: we always assume that each time a train service requires a resource $s\in S$ it indeed exploits one unit of the capacity $u_s$;
	\item a set $M$ of operations for trains. For each $m \in M$, we are given:
	\begin{itemize}
		\item a set of resources $S_{m} \subseteq S$ that are exploited by that operation;
		\item a completion time $\lambda_{m} \geq 0$;
		\item the maximum time $\gamma_{m} \geq 0$ that a train is allowed to wait after the completion of $m$ and before the start of the next operation: we let $\gamma_{m}:=\infty$ if this time is unbounded.
		\end{itemize}
	The set $M$ includes an operation $\alpha$ and an operation $\delta$ that respectively correspond to the arrival and the departure from the yard. We have: $S_{\alpha} = S_{\delta} = \emptyset$ and $\lambda_{\alpha} = \gamma_{\alpha} = \lambda_{\delta} = \gamma_{\delta} = 0$;
    \item a set $T$ and a set $T'$ of trains services, where $T$ represent current trains services and $T'$ represent new potential trains services. For each $t \in T \cup T'$ we are given:
	\begin{itemize}
        \item a set of plans $P(t)$, where each plan $p \in P(t)$ is an ordered sequence of operations: i.e., $p = (p(1), p(2), ..., p(n_{p}))$, with $p(i) \in M$ for each $i\in [n_{p}]$. It is always the case that $p(1)=\alpha$ and $p(n_{p})=\delta$.
	\item a time window $[\underline{a}_t, \overline{a}_t]$ for the arrival time and a time window $[\underline{q}_t, \overline{q}_t]$ for the departure time: note that for each $t\in T$ $\underline{a}_t =\overline{a}_t$ and $\underline{q}_t = \overline{q}_t$, i.e., the time windows reduce respectively to the arrival and departure time in the current timetable.
	\end{itemize}
\end{itemize}

\begin{definition}\label{def:schedule}A \textit{schedule} for a train service $t \in T\cup T'$ is a pair $(p, \sigma)$, where $p= (p(1), p(2), ..., p(n_{p}))$ is a plan in $P(t)$ and $\sigma: p \mapsto [\underline{a}_t, \overline{q}_t]$ assigns to each operation $p(i)$ a starting time  $\sigma(p(i))$. 

The schedule $(p, \sigma)$ is feasible if $\sigma$ is consistent with both the arrival and departure times of $t$ and the completion of each operation, i.e.:
\begin{itemize}
	\item $\sigma(p(1)) \in [\underline{a}_t, \overline{a}_t]$ {\em and} $\sigma(p(n_p)) \in [\underline{q}_t, \overline{q}_t]$;
	\item $\sigma(p(i+1)) \geq \sigma(p(i)) + \lambda_{p(i)}$ {\em and}  $\sigma(p(i+1)) \leq \sigma(p(i)) + \lambda_{p(i)} + \gamma_{p(i)}$, $i = 1, \ldots n_{p} - 1$.
\end{itemize}
\end{definition}

Definition~\ref{def:schedule} deals with some consistency properties of the schedule of a {\em single} train service. The next definition concerns some consistency properties for {\em sets} of train services. The rationale is that {\em different} trains services may require a {\em same} resource $s\in S$ and we must therefore ensure that the capacity $u_s$ of that resource is not violated.

\begin{definition}\label{feasi_sched}
Given a set of train services $\overline T$, a set of schedules $\mathcal{P} = \{(p^t, \sigma^t), t \in \overline T\}$ is {\em feasible} if each schedule $(p^t, \sigma^t)$ is feasible and the schedules are such that the capacity of no resource is violated at any time. Namely, if for a resource $s \in S$ and time $\theta \in [\min_{t\in \overline T}\underline{a}_t, \max_{t\in \overline T}\overline{q}_t]$, we let $\overline T|_{\theta, s}^{\mathcal{P}}$ denote the set of train services using $s$ at time $\theta$, i.e., 
$$\overline T|_{\theta, s}^{\mathcal{P}} = \{t \in \overline T: s \in S_{p^t(i)} \mbox{\ and\ } \theta\in [\sigma^t(p^t(i)), \sigma^t(p^t(i+1)), \mbox{\ for\ some\ }  i\in [n_{p^t} - 1]\},$$
then the following must hold: $|\ \overline T|_{\theta, s}^{\mathcal{P}}\ | \leq u_s$.
\end{definition}

\medskip\noindent
\textbf{The Yard Saturation Problem} ({\sc ysp}) Given $(S, M, T, T')$ defined as above, find a set of schedules $\mathcal{P}= \{(p^t, \sigma^t), t \in T\cup T^*\}$ that is feasible, with $T^*\subseteq T'$ and $|T^*|$ maximum.

\smallskip
It is possibly the case that some additional constraints are defined for particular instances of the {\sc ysp} on top of those implicit in its definition. E.g., for the Marzaglia transshipment yard, there are constraints due to rresource unavailability (see Section~\ref{section:unavailable}) and average resource consumption (see Section~\ref{section:robustness}). While these constraints impact on the computational performances of the optimization model, they do not alter too much the general structure of the {\sc ysp}, so we postpone their definition. To the contrary, we devote the next section to discussing a key variation on the {\sc ysp}.

\subsection{Periodic Yard Saturation Problem}
\label{section:periodic_scheduling}

Timetables and schedules that are periodic are ideal to keep the daily management of a yard easier. Therefore timetables and schedules at a yard  {\em are} periodic and in order to evaluate the practical capacity we must take into account this issue.
In this section, we therefore add periodic constraints to the {\sc ysp}. In particular, we assume that we are given a period $\tau$ that we want to enforce so that each train service $t$ presents itself every $\tau$ time instants and has a schedule that is periodic with period $\tau$. We will therefore deal with multiple copy (trains) of each train service, that it is convenient to define as {\em replicas} of a first {\em base} replica defined for the period $[0, \tau)$.

An instance of the \textit{Periodic Yard Saturation Problem} ({\sc pysp}) is  defined by a tuple $(S, M, T, T', \tau)$, and it is convenient to think of it as the replica of the $(S, M, T, T')$ instance of the {\sc ysp} over time with period $\tau$. Therefore, while e.g. the instance $(S, M, T, T')$ provides for some train service $t \in T\cup T'$ a time window $[\underline{a}_t, \overline{a}_t)  \subseteq [0, \tau)$ for the arrival time,  that is intended for the first replica of $t$ and we should infer analogous bounds for the next replicas of $t$: e.g. in this case, the time window for the second replica is $[\underline{a}_t +\tau, \overline{a}_t + \tau)  \subseteq [\tau, 2\tau)$ etc. 

We extend this idea to schedules. We will then compute a schedule for the first replica of each train service $t$, as in Definition~\ref{feasi_sched}, but we will then assume that the same schedule, shifted by $h\cdot \tau$, holds for each replica of $t$. This leads to the following definition, that generalizes Definition~\ref{feasi_sched}, where we characterize when a periodic schedule is {\em feasible}:

\begin{definition}\label{feasi_sched_period}
Given a subset of train services $\overline T$, a set of schedules $\mathcal{P} = \{(p^t, \sigma^t), t \in \overline T\}$ is feasible and  periodic of period $\tau$ if each schedule $(p^t, \sigma^t)$ is feasible and the schedules are such that the capacity of no resource is violated at any time, even if we take into account replicas of each schedule in $\mathcal{P}$.
Namely, for a resource $s \in S$ and time $\theta \geq 0$, let:
\begin{align*}
\overline T|_{\theta, s, \tau}^{\mathcal{P}} = \{ & (t, h),  t\in \overline T, h \in \mathbb{N}_{\geq 0} : \exists ~ i \in [n_{p^t} - 1] s.t. ~ s \in S_{p^t(i)} \\
 & \mbox{\ and\ } \theta\in [\sigma^t(p^t(i)) + h \tau , \sigma^t(p^t(i+1)) + h \tau)\}
\end{align*}
the set of replicas of train services using $s$ at time $\theta$.
Then the following must hold: $$|\overline T|_{\theta, s, \tau}^{\mathcal{P}}\ | \leq u_s.$$
\end{definition}

\medskip\noindent
\textbf{The Periodic Yard Saturation Problem} ({\sc pysp}) Given $(S, M, T, T', \tau)$ defined as above, find a set of schedules $\mathcal{P}= \{(p^t, \sigma^t), t \in T\cup T^*\}$ that is feasible and periodic, with $T^*\subseteq T'$ and $|T^*|$ maximum.

\smallskip
The problem seems quite challenging since we are now dealing with an unbounded time horizon. However we may restrict to a finite number of periods and replicas of a same train service, and in practice reduce to the solution of a ``time-extended" (non-periodic) Yard Saturation Problem. This goes as follows. First, let us define the {\em makespan} $\zeta(T\cup T')$: $$\zeta(T\cup T') = \max_{t \in T\cup T'} \overline{q}_t - \min_{t \in T\cup T'} \underline{a}_t$$

Suppose that $\zeta(T\cup T') \leq \tau$. It is then possible to obtain a solution to an instance of the {\sc pysp} by simply replicating a solution to the corresponding instance of {\sc ysp}. However, if $\zeta(T\cup T') > \tau$, then it is possible that two replicas from different periods compete for a same resource. In this case, let $k(T\cup T'):= \lceil \frac{\zeta(T\cup T')}{\tau} \rceil$. Now, following the definition of $k(T\cup T')$, if the schedule of the $i$-th replica of a train service $t$ and that of the $j$-th replica of a train service $t'$ overlap, and  possibly compete for a resource, it must be the case that $|i-j|\leq k(T\cup T') -1 $. Suppose without loss of generality that $i<j$. Then the conflict between those replicas has to be solved as the conflict between the first replica of $t$ and the $(j-(i-1))-th$ replica of $t'$. 

It follows therefore that we may restrict to an instance of the {\sc pysp} by simply replicating a solution to an instance of {\sc ysp} that is restricted to:
\begin{itemize}
    \item $k(T \cup T')$ replicas of each train service $t\in T\cup T$;
    \item the conflicts between the first replica of $t$ and the $j$-th replica of $t'$, for any $t$ and $t'\in T\cup T$ and any $j\in [k(T \cup T']$.
\end{itemize} 

For the sake of clarity, by ``replicating a solution" we mean that:

\begin{itemize}
    \item the schedule for each replica of $t$ is the same -- shifted by $h\cdot \tau$ -- as the schedule for the first replica of $t$;
    \item the plan for each replica of $t$ is the same  as the plan for the first replica of $t$;
    \item any conflict between the $i$-th replica of $t$ and the $j$-th replica of $t'$ with $i>j$ has to be solved as $(t,1)$ and $(t',j-(i-1))$.
\end{itemize}

\section{A disjunctive graph for the Yard Saturation Problem}
\label{section:disjunctive_graphs}

In this section, we define a suitable disjunctive graph for the solution of the Yard Saturation Problem. Disjunctive graphs, a popular tool for modeling job-shop scheduling problems, were first introduced in \cite{roy_1964}, see also \cite{balas1}, \cite{balas2}. In the context of railway transportation disjunctive graphs were popularized in \cite{mascis_2002} where they were called alternative graphs.

In a disjunctive graph nodes, representing tasks or operations to be performed, are connected by directed arcs representing either precedence constraints or timing constraints between operations. We in particular assume that an arc $(p, q)$ with length $l(p,q)$ models the constraints: \begin{equation}
\label{scheduling arcs}
t_q \geq t_p + l(p,q)
\end{equation}
where $t_p$ and $t_q$ are respectively the starting time of operation $p$ and $q$ (and $l(p,q)$  is e.g. the duration of $p$). We point out that we allow for some arc $(j,i)$ negative length, i.e., $l(j,i)<0$. This is useful e.g. because if there is also the arc $(i, j)$ with length $l(i,j)>0$, then the arcs $(i, j)$ and $(j,i)$ model the constraints $l(i,j) \leq t_j - t_i \leq |l(j,i)|$, which allow to model time windows constraints.

We deal with two types of precedence constraints and therefore arcs: 
\begin{itemize}
    \item[(i)] {\em strict precedence constraint}, this is the case of a task $p$ that is the predecessor of a task $q$ in some operation plan. It is modeled by an arc $(p, q)$ usually with a length $l(p, q)>0$ representing the duration of $p$;
    \item[(ii)] {\em disjunctive precedence constraint}, this is the case of two operations $p$ and $q$ that are independent and can be performed in either order, but not simultaneously e.g. because they both require a same resource. It is modeled by a pair of suitable {\em disjunctive} arcs, i.e., every feasible solution will select exactly one of them (see below for details).
\end{itemize} 

We point out that constraint $(ii)$ is slightly different when the operations $p$ and $q$ compete for the use of a same resource that is available in multiple units, i.e., it has capacity larger than one. In this case, first considered in \cite{lamorgese_2015}, it is possible that $p$ and $q$ are performed simultaneously (possibly for some short period of time) and this is now represented by four capacity arcs: two {\em precedence} arcs are defined as above in $(ii)$, two {\em meeting} arcs model the fact that $p$ and $q$ overlap, possibly for just a short period of time. When we solve the disjunction we will either select exactly one precedence arc or both the meeting arcs: see below for more details. Note that we must however guarantee, via a suitable additional constraints (in our case this will be Inequality~\eqref{eq:capacity_constraints}) in Section~\ref{section:rfi_milp}, that the number of operations that are simultaneously using the resource is not larger than its capacity.

\subsection{Single plan for each train service}\label{spets}

We are now ready to associate to an instance $(S, T, T', M)$ of the (non-periodic) {\sc yps}, a suitable disjunctive graph $G(V, A, l)$, with $l:A\rightarrow{\cal Z}_+$. For the sake of simplicity, we initially assume that each train service has just one possible plan. Therefore, for each train $t \in T\cup T'$, $P(t) = \{p^t\}$ and $p^t = (p^t(1), p^t(2), ..., p^t(n_{p^t}))$. 

\paragraph{Nodes}
We associate with each operation $p^t(j)$ in the plan $p^t$ of a train service $t\in T\cup T'$ a node $v \equiv p^t(j)$ -- i.e., $p^t(j)$ denotes both the operation and the node, and with a little abuse of notation we will also let $S_v = S_{p^t(j)}$.  We then define a last additional node $o$ that stands for the beginning of the time horizon. 

$$ V = \{o\} \cup \bigcup_{t \in T \cup T'} ~ V(t), \text{\ with\ } V(t): = \{p^t(1), p^t(2), ..., p^t(n_{p^t})\}.$$

\paragraph{Arcs}
The set of arcs $A$ is such that $A = A_s \cup A_{c}$, where the set $A_s$ models strict precedence constraints and the set $A_{c}$ models disjunctive precedence constraints. 

\smallskip\noindent
{\bf The set $A_s$} The set $A_s$ is such that $A_s =\bigcup_{t \in T\cup T'} ~ A(t)$ with $ A(t) = \overrightarrow A(t) \cup \overleftarrow A(t) \cup A(t)_{a} \cup A(t)_{d}$:
\begin{itemize}
\item $\overrightarrow A(t) = \{(p^t(j), p^t(j+1)), j \in [n_{p^t}-1]\}$, with $l((p^t(j), p^t(j+1))): = \lambda_{p^t(j)}$: these arcs impose precedence between consecutive operations; 
\item $\overleftarrow A(t) = \{(p^t(j+1), p^t(j)), j \in [n_{p^t}-1]\}$, with $l((p^t(j+1), p^t(j))): = - (\lambda_{p^t(j)} + \gamma_{p^t(j)})$: these arcs model constraints due to the maximum waiting time after each operation;
\item $A(t)_{a} = \{(o, p^t(1)), (p^t(1), o)\}$, having length respectively $\underline{a}^t$ and  $-\overline{a}^t$: these arcs model constraints on the arrival time of a train service;
\item $A(t)_{d} = \{(o, p^t(n_{p^t})), (p^t(n_{p^t}), o)\}$, having length respectively $\underline{q}^t$ and  $-\overline{q}^t$: these arcs model constraints on the departure time of a train service.
\end{itemize}

\begin{figure}
\begin{center}
   \begin{tikzpicture}[main/.style = {draw, circle, minimum size=1.2cm}] 
\pgfmathsetmacro{\l}{2}
\node[main] (1) at (1*\l,-1*\l) {$o$};
\node[main] (2) at (2.5*\l,-1*\l) {$p^t(1)$};
\node[main] (3) at (4*\l,-1*\l) {$p^t(2)$};
\node[main] (4) at (5.5*\l,-1*\l) {$p^t(3)$};
\node[main] (5) at (7*\l,-1*\l) {$p^t(4)$};
\node[main] (6) at (8.5*\l,-1*\l) {$p^t(5)$};

\begin{scope}[every edge/.style={draw=black, thick}]
\path [-{>[scale=\arrowHeadScale]}] (2) edge node [pos=0.5, sloped, above]{$\lambda_{p^t(1)}$} (3);
\path [-{>[scale=\arrowHeadScale]}] (3) edge node [pos=0.5, sloped, above]{$\lambda_{p^t(2)}$} (4);
\path [-{>[scale=\arrowHeadScale]}] (4) edge [bend left] node [pos=0.5, sloped, above]{$\lambda_{p^t(3)}$} (5);
\path [-{>[scale=\arrowHeadScale]}] (5) edge node [pos=0.5, above]{$\lambda_{p^t(4)}$} (6);
\end{scope}

\begin{scope}[every edge/.style={draw=green, thick}]
\path [-{>[scale=\arrowHeadScale]}] (1) edge [bend left] node [pos=0.5, sloped, below]{$\underline{a}^t$} (2);
\path [-{>[scale=\arrowHeadScale]}] (2) edge [bend left] node [pos=0.5, sloped, above]{$-\overline{a}^t$} (1);
\end{scope}

\begin{scope}[every edge/.style={draw=purple, thick}]
\path [-{>[scale=\arrowHeadScale]}] (1.north) edge [bend left] node [pos=0.5, sloped, above]{$\underline{q}^t$} (6.north);
\path [-{>[scale=\arrowHeadScale]}] (6.south) edge [bend left] node [pos=0.5, sloped, below]{$-\overline{q}^t$} (1.south);
\end{scope}

\begin{scope}[every edge/.style={draw=cyan, thick}]
\path [-{>[scale=\arrowHeadScale]}] (5) edge [bend left] node [pos=0.5, sloped, below]{$-\lambda_{p^t(3)}-\gamma_{p^t(2)}$} (4);
\end{scope}

\end{tikzpicture}  
\end{center}
\caption{The set of arcs $A(t)$ for a train service $t$ with $p^t = (p^t(1), \ldots, p^t(5))$. The arcs in $\protect\overrightarrow{A(t)}$ are in black, those in $\protect\overleftarrow{A(t)}$ in cyan, those in $A(t)_{a}$ in green, and finally those in $A(t)_{d}$ in purple. Note that $\gamma_{p^t(1)} = \gamma_{p^t(2)} = \gamma_{p^t(4)} =\infty$, so we are omitting the corresponding arcs in $\protect\overleftarrow{A(t)}$.}
\label{fig:arcs}
\end{figure}
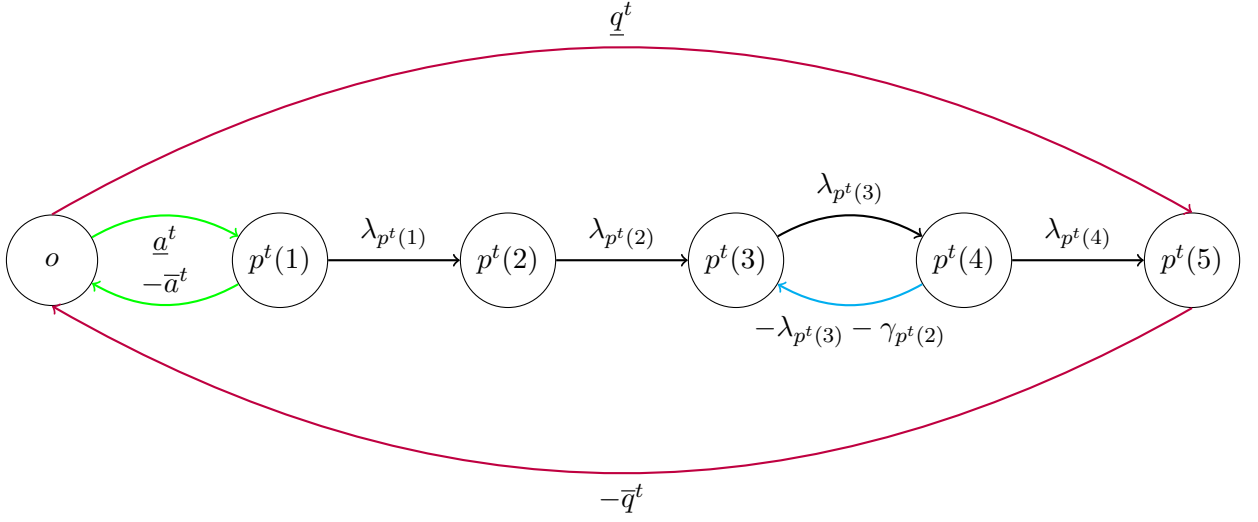

Figure~\ref{fig:arcs} shows how set $A(t)$ might look like for a single train. 

\smallskip\noindent
{\bf The set $A_c$} The set of arcs $A_{c}$ deals with potential conflicts between operations in the plan of different train services that require a same resource. Its definition requires a few tools. First, we say that two operations $m$ and $m'\in M$ are in \textit{conflict} if $S_m\cap S_{m'}\neq \emptyset$. 

\begin{definition}\label{setF}
We denote by $F$ the set of (unordered) pairs of nodes from $V(t)$ that correspond to operations that are in conflict. Namely: 
$$F = \{\{u, v\}: u\in V(t), v\in V(t'), t, t' \in T\cup T', S_{u} \cap S_{v} \neq \emptyset\}.$$
\end{definition}

Note that it is possibly the case that $\{u, v\}\in F$, with both $u$ and $v$ in the same set $V(t)$ for some $t\in T\cup T'$ (i.e., $t=t'$ in the above definition). In this case, the precedence between $u$ and $v$ is also strict, however defining $F$ as so is convenient for the problem in the periodic setting (see Section~\ref{periodic}) where there are multiple replicas of a same train service and we have to decide the precedence between operations corresponding to $u$ and $v$ from different replicas.

\begin{definition}\label{conflict_arcs}
Let $\{u, v\}\in F$ with $u\equiv p^t(i)$ and $v\equiv p^{t'}(j)$. We associate with $\{u, v\}$ the set $A(\{u, v\})$ that is made of the following four arcs:
\begin{itemize}
\item $(p^t(i+1), p^{t'}(j))$ and $(p^{t'}(j+1), p^t(i))$, with length $\epsilon$: these are called \textit{precedence arcs} and respectively models that  $p^t(i)$ (resp. $p^{t'}(j)$) is completed before $p^{t'}(j)$ (resp. $p^t(i)$) starts;
\item $(p^t(i), p^{t'}(j+1))$ and $(p^{t'}(j), p^t(i+1))$, with length $-\epsilon$: these are called \textit{meeting arcs} and model that $p^t(i)$ and $p^{t'}(j)$ may overlap,
\end{itemize}
\end{definition}
where $\epsilon$ is a small positive value.

\begin{figure}
\begin{center}
\begin{tikzpicture}[main/.style = {draw, circle, minimum size=1.3cm}] 
\pgfmathsetmacro{\l}{2}
\node[main] (1) at (0*\l,0*\l) {$o$};
\node[main] (2) at (1*\l,1*\l) {$p^{t}(1)$};
\node[main] (3) at (2*\l,1*\l) {$p^{t}(2)$};
\node[main] (4) at (3*\l,1*\l) {$p^{t}(3)$};
\node[main] (5) at (4*\l,1*\l) {$p^{t}(4)$};
\node[main] (6) at (1*\l,-1*\l) {$p^{t'}(1)$};
\node[main] (7) at (2*\l,-1*\l) {$p^{t'}(2)$};
\node[main] (8) at (3*\l,-1*\l) {$p^{t'}(3)$};
\node[main] (9) at (4*\l,-1*\l) {$p^{t'}(4)$};

\begin{scope}[every edge/.style={draw=red}]
\path [->] (4) edge node [red, pos=0.5, sloped, above]{$\epsilon$} (7);
\end{scope}
\begin{scope}[every edge/.style={draw=blue}]
\path [->] (8) edge node [blue, pos=0.5, sloped, above]{$\epsilon$} (3);
\end{scope}
\begin{scope}[every edge/.style={draw=orange}]
\path [->] (3) edge [bend left] node [orange, pos=0.5, sloped, above]{$-\epsilon$} (8);
\path [->] (7) edge [bend left] node [orange, pos=0.5, sloped, above]{$-\epsilon$} (4);
\end{scope}
\begin{scope}[every edge/.style={draw=black}]
\path [->] (1) edge node [pos=0.5, sloped, above]{} (2);
\path [->] (2) edge node [pos=0.5, sloped, above]{} (3);
\path [->] (3) edge node [pos=0.5, sloped, above]{} (4);
\path [->] (4) edge node [pos=0.5, sloped, above]{} (5);
\path [->] (1) edge node [pos=0.5, sloped, above]{} (6);
\path [->] (6) edge node [pos=0.5, sloped, above]{} (7);
\path [->] (7) edge node [pos=0.5, sloped, above]{} (8);
\path [->] (8) edge node [pos=0.5, sloped, above]{} (9);
\end{scope}
\end{tikzpicture} 
\end{center}
\caption{The operations corresponding to $p^t(2)$ and $p^{t'}(2)$ are in conflict. Precedence arcs are red and blue; meeting arcs orange. For the sake of simplicity, some arcs in $A(t)$ and $A(t')$ are omitted.} 
\label{fig:conflict}
\end{figure}
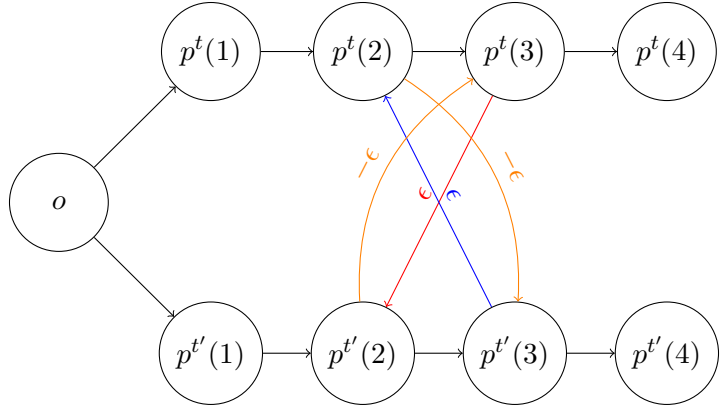

These arcs in $A(\{u, v\})$ completely describe the possible interactions between the conflicting operations $u$ and $v$. The precedence arcs model the case where either operation is completed before the other starts. The meeting arcs model the case where the operations are performed simultaneously, which is possible when the resources in $S(p^t(i)) \cap S(p^{t'}(j))$ have capacity greater than one. See Figure~\ref{fig:conflict} for an example.

Note that assigning a non-zero length to precedence and meeting arcs ensures that once a resource is freed it can not be \textit{immediately} occupied by another train, enforcing the so-called \textit{no-swap constraints} \cite{mascis_2002}.

The set of arcs $A_c$ is therefore the following:
\begin{equation}
\label{conflict arcs}
A_c = \bigcup_{\{u, v\}\in F} A(\{u, v\}).
\end{equation}.
and in any solution to the {\sc ysp} for each $\{u, v\}\in F$ we will select  {\em either one of the precedence arcs or both the meeting arcs} in $A(\{u, v\})$ (see Section~\ref{section:rfi_milp} for more details).

\subsection{A disjunctive graph for the Periodic Yard Saturation Problem}\label{periodic}
We now discuss how the disjunctive graph can be extended in order to deal with an instance $(S, M, T, T', \tau)$ of the Periodic Yard Saturation Problem. We will indeed only sketch how it is possible to generalize the disjunctive graph presented in the previous section to the periodic setting, a more detailed discussion can be found in \ref{appendix:periodic_scheduling}. We still assume that for each train service there is just one possible plan. 

Following the discussion in Section~\ref{section:periodic_scheduling}, we need to take into account $k$ replicas of each train service, with $k:=  k(T\cup T')$. However, our choice is that of dealing with the same set of nodes $V$ that was defined in the previous section, and the same will happen for the set of arcs $A_s$. We will indeed deal with the different replicas through the arcs of the set $A_c$: recall that in a disjunctive model constraints correspond to arcs and we may therefore impose all constraints we need between replicas of $u$ and $v$ adding more arcs, each with a suitable length, to the set $A_c$. This goes as follows. The set $A_c$, it is still defined as in \eqref{conflict arcs}, however each set $A(\{u, v\})$ is now made of $4(2k-1)$ arcs. Following again the discussion in Section~\ref{section:periodic_scheduling}, for each pair of conflicting nodes $\{u, v\}\in F$, with $u\in V(t)$ and $v\in V(t')$, we need to decide the precedences between the operation $u$ in the 
first replica of $t$, that we simply denote as $u$, and the operation $v$ in each replica of $t'$, that we denote as $v, v^2, \ldots, v^k$, as well as the precedences between $v$ and the operation $u$ in each replica of $t$, namely $u, u^2, \ldots, u^k$. So, for each pair of conflicting nodes $\{u, v\}\in F$, we need to decide $2k-1$ precedences and each one of this precedence can be modeled through the 4 arcs that have been introduced in the previous section. Altogether, there will be $4(2k-1)$ arcs, however the length of these arcs will depend on which replica (either of $u$ or $v$) is involved. A more detailed discussion can be found in \ref{appendix:periodic_scheduling}.

\subsection{Multiple plans for each train service}
\label{section:multiple_plans}

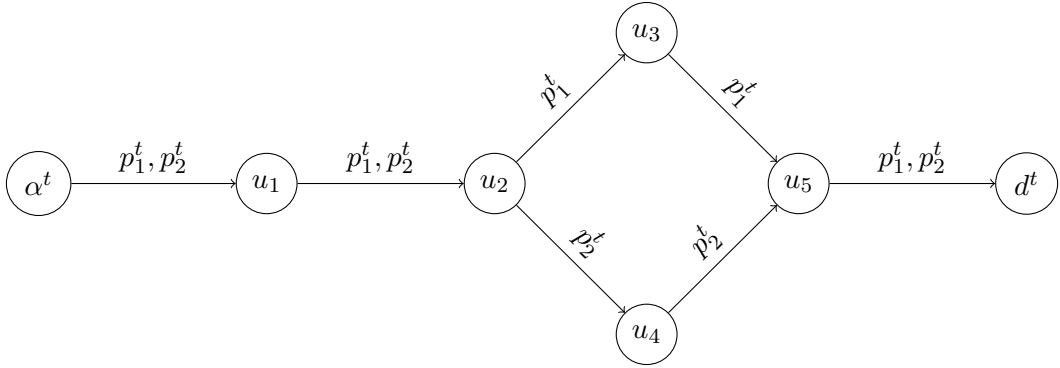
\begin{figure}
\begin{center}
\begin{tikzpicture}[main/.style = {draw, circle}] 
\pgfmathsetmacro{\l}{2}
\node[main] (8) at (0*\l,-1*\l) {$\alpha^t$};
\node[main] (2) at (1.5*\l,-1*\l) {$u_1$};
\node[main] (3) at (3*\l,-1*\l) {$u_2$};
\node[main] (4) at (4*\l,0*\l) {$u_3$};
\node[main] (5) at (4*\l,-2*\l) {$u_4$};
\node[main] (6) at (5*\l,-1*\l) {$u_5$};
\node[main] (7) at (6.5*\l, -1*\l) {$d^t$};

\begin{scope}[every edge/.style={draw=black}]
\path [->] (8) edge node [pos=0.5, sloped, above]{$p^t_1, p^t_2$} (2);
\path [->] (2) edge node [pos=0.5, sloped, above]{$p^t_1, p^t_2$} (3);
\path [->] (3) edge node [pos=0.5, sloped, above]{$p^t_1$} (4);
\path [->] (3) edge node [pos=0.5, sloped, above]{$p^t_2$} (5);
\path [->] (4) edge node [pos=0.5, sloped, above]{$p^t_1$} (6);
\path [->] (5) edge node [pos=0.5, sloped, above]{$p^t_2$} (6);
\path [->] (6) edge node [pos=0.5, above]{$p^t_1, p^t_2$} (7);
\end{scope}
\end{tikzpicture}
\end{center}
\caption{Taking advantage of common operations in two different plans $p^t_1$ and $p^t_2$ for a same train service. With some abuse of notation, we have that $p^t_1 = (\alpha^t, u_1, u_2, u_3, u_5, \delta^t)$, while $p^t_2 = (\alpha^t, u_1, u_2, u_4, u_5, \delta^t)$. For simplicity, in this picture we draw only the arcs in $\protect\overrightarrow{A(t)}$.}
\label{fig:multiple_plans}
\end{figure}

We now discuss how to relax the hypothesis $|P(t)| = 1 ~ \forall t \in T$. A simple strategy would be that of replacing the single plan $p(t)$ considered so far with a set of parallel plans for each train.
However, that is not the most convenient choice. Indeed, as we already discussed in Section~\ref{sc:trans}, different plans of a same train service usually share several common operations; therefore, replicating every operation for each feasible plan determines an unnecessary growth of the size of the disjunctive graph and introduces symmetries. We therefore start from the latter representation but ``shrink" nodes and arcs so as to take advantage of that feature.
Figure~\ref{fig:multiple_plans} provides an intuition for this, note that both the set of nodes and the set of arcs need to be suitably defined and they are larger than those defined for the single plan case. Additional details are provided in \ref{appendix:multiple_plans}.

\section{A Mixed Integer Linear Programming formulation for the Periodic Yard Saturation Problem}
\label{section:rfi_milp}

In this section, we present a Mixed Integer Linear Programming formulation for the {\sc pysp}. This formulation builds upon the disjunctive graph that was presented in Section~\ref{periodic}. Our formulation will be a {\em Big-M formulation}, where we use a suitable large constant $M$ to linearize the disjunctive constraint arising from disjunctive arcs.

\smallskip
We start with introducing the main constraints and variables of our model. Note that while our discussion in Section~\ref{periodic} was the case of single plan for each train, in the following we implicitly deal with the case where there are more plans for each train and refer to the \ref{appendix:unavail} in the few cases where we need particular tools.

\paragraph{Plan selection}
We must choose one plan for each train service $t$ that is served. We therefore define:
\begin{itemize}
	\item for each train service $t \in T\cup T'$, a variable $\phi_t \in \{0, 1\}$ that is set to 1 if and only if the train service $t$ is active. Note that: 
\begin{equation}
\label{eq:fixed_trains}
\phi_t = 1 \qquad t \in T
\end{equation}
	\item for each train service $t \in T\cup T'$ and each plan $p \in P(t)$, a variable $w_{p} \in \{0, 1\}$ that is set to 1 if only if the train service $t$ is active and the plan $p$ is selected for $t$.
\end{itemize}
The following constraint connect the above variables: 
\begin{equation}
\label{eq:plan_constraints}
\sum_{p \in P(t)} w_{p} = \phi_t \qquad t \in T\cup T'.
\end{equation}

\paragraph{Active nodes}
We introduce the notion of an \textit{active node}: a node $u \in V(t)$ is active if it belongs to an active plan; the origin node $o$ is always active. Then we associate a binary variable $x_u \in \{0, 1\}$ for each node $u \in V$ and: 
\begin{equation}
\label{eq:node_activation}
x_u = \begin{cases}
        1 \qquad \qquad \quad ~ \text{if $u = o$,}\\
        \sum_{p \in P(t): u\in V(p)} w_{p} \qquad \text{if\ } u\in V(t), t \in T\cup T'
        \end{cases} \qquad u \in V
\end{equation}
where $V(p)$ is the set of nodes associated with plan $p$.

\paragraph{Active arcs}
We also introduce the notion of an \textit{active arc}. Recall that in a disjunctive graph each arc induces a constraint (see Inequality \eqref{scheduling arcs}) but it is not always the case that that constraints is indeed active: e.g. because it corresponds to an arc that is part of a plan that has not been selected, or because it corresponds to a an arc for a disjunction that has been solved as to remove that constraint.

\smallskip\noindent
{\bf Active arcs in $A_s$}. We first deal with the arcs in $A_s$ and in particular with arcs in $A(t)$, for each $t\in T\cup T'$. Recall that $A(t) = \overrightarrow A(t) \cup \overleftarrow A(t) \cup A(t)_{a} \cup A(t)_{d}$:, we have that $a \in A(t)_{a} \cup A(t)_{d}$ is active if and only if $t$ is served by the yard, while $a \in \overrightarrow{A(t)} \cup \overleftarrow{A(t)}$ is active if and only if $t$ is served and $a$ belongs to a plan for $t$ that is active. Namely:
\begin{equation}
x_a = \begin{cases}
        \phi_t  ~~ \text{if $a \in A(t)_{a} \cup A(t)_{d}$,}\\
         \sum_{p \in P(t): a\in A(p)} w_{p} ~~ \text{otherwise}.
        \end{cases} \qquad a \in A(t), t \in T \cup T'
\end{equation}
Where $A(p)$ is the set of arcs associated with plan $p$.

\smallskip\noindent
{\bf Active arcs in $A_c$}.
Recall that for each pair $\{u, v\}\in F$ we define $2k-1$ tuples each made of four arcs each, with each tuple modeling either the precedences between $u$ and some replica $v^i$ of $v$ or the precedences between $v$ and some replica $u^j$ of $u$, with both $i$ and $j= 1..k$. We therefore define, e.g. for each pair $(u, v^i)$, three binary variables, respectively $y_{uv^i}$, $y_{v^iu}$ and $z_{\{uv^i\}}$ since the meeting arcs are either both active or both not active.
The following constraint connect the above with $x_{u}$ and $x_{v}$, for each $i=1..k$:
\begin{equation}
\label{eq:selection_constraints}
\begin{split}
    y_{uv^i} + y_{v^iu} + z_{\{uv^i\}} \geq x_u + x_v - 1 \\
    y_{vu^i} + y_{u^iv} + z_{\{vu^i\}} \geq x_u + x_v - 1
\end{split}
\qquad \{u, v\}\in F, i= 1..k
\end{equation}

 We now associate an activation variable $x_a$ also with each conflict arc $a\in A_c$. When there is a single plan for each train service, the $y$ and $z$-variables defined above will do the job. In other words, it will be enough to set the value of $x_a$ equal either to the value of a suitable $y$-variable, if $a$ is a precedence arc, or to the value of a suitable $z$-variable if a is meeting arc. However, when there are more plans for a same train service, things get more complicated. We therefore postpone the formal statement of these constraints to \ref{appendix:multiple_plans_milp}. 

\paragraph{Scheduling variables and constraints}
We associate with each node $v\in V$ a continuous variable $\sigma(v) \geq 0$ representing the starting time of the operation associated with node $v$. We set $\sigma(o) = 0$ and then write the following standard precedence constraint:
\begin{equation}
\label{precedence_constraint}
\sigma(v) - \sigma(u) \geq l_{a} - M(1 - x_{a}) \qquad a\in A, v\ \text{tail of }a, u\ \text{head of }a.
\end{equation}

\paragraph{Capacity constraints}
The $z$-variables allow for the possibility of multiple operations using the same resource at the same time. We must check that this is done without exceeding the capacity of the resource.
We denote by $O_s$ is the set of all operations (and replicas) requiring a resource $s\in S$: $$O_s = \{u^i, i = 1, ..., k, u \in p^t, p^t \in P(t), t \in T\cup T': s \in S_u\}.$$
According to the next constraint, each resource $s$ is used by at most $u_s$ operations (and replicas) at a time:
\begin{equation}
\label{eq:capacity_constraints}
\sum_{\{u, v^i\} \subseteq Q} z_{\{uv^i\}} \leq \binom{u_s + 1}{2} - 1 \qquad \forall s \in S, \forall Q \subseteq O_s: |Q| = u_s + 1
\end{equation}

\paragraph{Objective function}
Finally, the objective function that we want to maximize is the number of served train, i.e. $\sum_{t \in T\cup T'} \phi_t$.

\subsection{Capacity cuts}
\label{section:analytical_bounds}

Constraints~\eqref{eq:fixed_trains}-\eqref{eq:capacity_constraints} provide a first {\sc milp}-formulation for {\sc pysp}. We now show how it is possible to improve this formulation by some suitable {\em cutting planes}. By exploiting some information on the trains set, it is indeed possible to estimate an upper bound on the value $z^*$ of the optimal solution to an instance of the {\sc pysp}, also in presence of multiple plans and replicas.

Let $b$ be the total number of tracks in the transshipment yard and let $\underline{\rho}_t = \underline{q}_t - \overline{a}_t$ be the minimum stay of a train service $t \in T \cup T'$. Then a first upper bound for $z^*$ is given by: 

\begin{equation}
\label{eq:cap_risorse}
z^* \leq \lfloor \frac{b \tau }{min_{t \in T} \underline{\rho}_t} \rfloor
\end{equation}

Then, for each resource $s\in S$, let $\mu_s^{t, p}$ be the minimum occupation time of $s$ by $t$, if plan $p \in P(t)$ is the one selected by $t$. We then define $\mu_s^{t}$ as the minimum occupation time of $s$ by $t$, i.e. $\mu_s^{t} = min_{p \in P(t)} \mu_s^{t, p}$. Another upper bound for $z^*$ is then given by:
\begin{equation}
\label{eq:cap_risorsebis}
z^* \leq min_{s \in S} \lfloor \frac{u_s \tau}{min_{t \in T} \mu_s^t} \rfloor
\end{equation}

As we show in \ref{appendix:complete_milp}, these constraints can be easily linearized and then used as cutting planes for improving the quality of our formulation. In the Appendix we also provide the entire formulation for {\sc pysp} that takes into account all constraints discussed so far and a few other more technical constraints that we present in the following.

\subsection{Additional constraints}
In this section we discuss some additional constraints for the Yard Saturation Problem. These constraints are quite natural in the setting of the (Periodic) Yard Saturation Problem and in fact they were part of the model we developed for our case study in Marzaglia. In the following, we shortly discuss these constraints and refer to the Appendix for a formal statement of these constraints in terms of Mixed Integer Linear Programming.

\subsubsection{Resource unavailabilities}
\label{section:unavailable}
In a transshipment yard resources may be not available because in some time intervals because of maintenance or time-shifts. This is e.g. for our case study in Marzaglia, where the transshipment equipment is not always available because of work shifts of their operators. 

We therefore consider the (periodic) setting where some resource $s$ (possibly more than one) is not available for some regular interval times (e.g. in Marzaglia from 11 pm to 5 am). In \ref{appendix:unavail} we show that our model can be easily extended as to encode such constraints, also in presence of multiple plans and replicas.

\subsubsection{Bound on average utilization}
\label{section:robustness}
In order to increase the robustness of the solutions computed in this study, RFI suggested a simple rule of thumb: the average utilization of each resource should not be higher than $85\%$ of the full capacity. Once again our model can be easily extended as to encode such constraints, also in presence of multiple plans and replicas. This is again discussed in the \ref{appendix:robustness}.

\subsection{Solving the {\sc milp} model}
\label{section:algorithm_M}

The {\sc milp} formulation of the previous section can grow to a prohibitively large size, even for small instances. This is mostly due to Constraints \eqref{eq:capacity_constraints}. We therefore adopt a dynamic row and column generation iterative scheme in which at each iteration we consider only a subset of conflicting operations.
The resulting {\sc milp} problem has both a reduced set of constraints and a reduced set of variables.
Once we have solved this problem, the solution provides a set of schedules which in general is not feasible. To determine if it is feasible, we must solve a separation problem for each resource $s$. This can be done by constructing an interval graph for each resource $s$ and by checking whether there exists a clique of size $u_s + 1$ or more. Recall that on an interval graph it is possible to find the set of all maximal cliques in linear time \cite{rolf_1986}. If there are no such cliques, then the algorithm terminates with an optimal solution. Otherwise, we enlarge the considered subset of conflicting operations and continue with the next iteration. A more exhaustive discussion on the dynamic row and column generation scheme that we adopted can be found in \ref{appendix:algorithm}.

\smallskip
The first feasible solution found by the above scheme is of course optimal. However, when solving large instances, we could be interested in obtaining intermediate sub-optimal solutions. To address this issue we solve a sequence of feasibility problems where we ask for solution with larger and larger size: namely, as soon as we find a feasible solution that serves $f$ trains we add the constraint
$$\sum_{t \in T \cup T'} \phi_t \geq f + 1$$
The procedure terminates when either $f = |T \cup T'|$ or the current problem is unfeasible. In both cases, we conclude that $f$ is the optimum value.

\section{Computational Results for Marzaglia yard}
\label{section:rfi_results}

\begin{figure*}
        \centering
        \begin{subfigure}[b]{0.475\textwidth}
            \centering
            \includegraphics[width=\textwidth]{layout_01.png}
            \caption{Scenarios 0 and 1}
            \label{fig:layout_01}
        \end{subfigure}
        \vskip\baselineskip
        \begin{subfigure}[b]{0.475\textwidth}  
            \centering 
            \includegraphics[width=\textwidth]{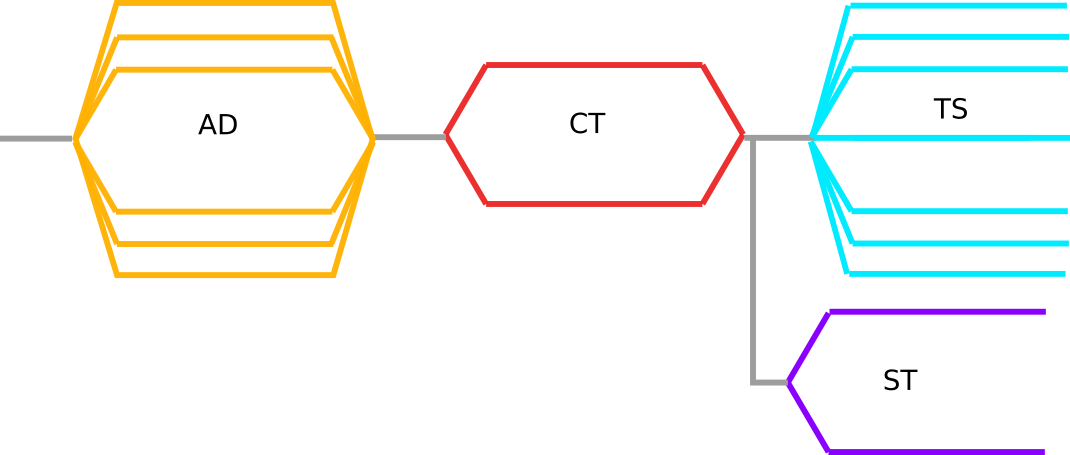}
            \caption{Scenario 2}
            \label{fig:layout_2}
        \end{subfigure}
        \hfill
        \begin{subfigure}[b]{0.475\textwidth}  
            \centering 
            \includegraphics[width=\textwidth]{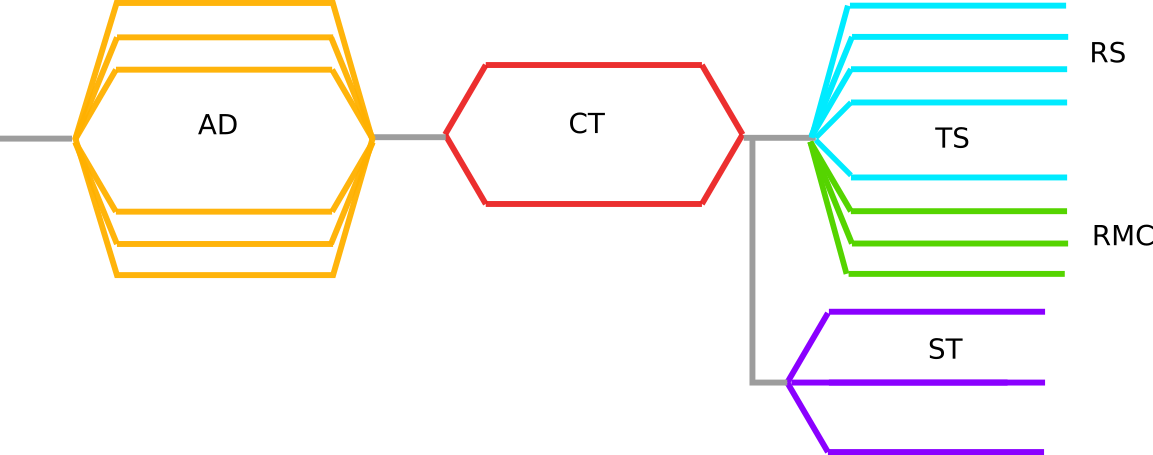}
            \caption{Scenario 3}
            \label{fig:layout_3}
        \end{subfigure}
        \caption{Infrastructural layout for the different scenarios.} 
        \label{fig:layouts}
    \end{figure*}

RFI was interested in evaluating the practical capacity of the transshipment yard of Marzaglia both in the current configuration and in three other configurations, corresponding to different investment plans. The current layout is depicted in Figure~\ref{fig:layouts} (a) and will be referred in the following as {\em Scenario} 0. In Scenario 0 loading and unloading operations are done by reach stackers that are operated from 5 am to 11 pm. The investment plans involve:
\begin{itemize}
    \item Scenario 1: the extension of the operation time for reach stackers for so that they can be operated 24 hours a day from Monday to Saturday;
    \item Scenario 2: the extension of the operation time for reach stackers as in Scenario 1 and the addition of some more parallel tracks to the current layout, see Figure~\ref{fig:layouts} (b); 
    \item Scenario 3: the extension of the operation time for reach stackers as in Scenario 1, the addition of some more parallel tracks to the current layout, see Figure~\ref{fig:layouts} (c), and the addition of a rail mounted crane spanning over three transshipment tracks.
\end{itemize}

The infrastructural layout associated with each scenario is represented in Figure~\ref{fig:layouts}. Across all scenarios, each train must move from the arrival/departure area to the transshipment area and back, and has the possibility of moving to the side track area either in its forward or backward journey, or both. This means that each train has 4 plans available. In addition, in Scenario 3, for each train there is also a choice on the transshipment equipment used, implying that for each train we have to choose among 8 different plans.

These scenarios are associated with increasing and significantly different costs. RFI was interested in understanding the impact of each investment plan on the practical capacity of the yard, with a special focus on two ambitious goals with respect to the current timetables for the arrivals, that is provided in Table~\ref{table:current_timetable}: 1) operating all current train services {\em daily} as to increase the number of trains a week from 47 to 66; 2) even more, adding 5 new train services to be operated daily -- and for which they suggested reasonable arrival and departure time windows -- as to increase the number of weekly trains to 96, therefore more than double the current number.

\paragraph{Scenarios 0 and 1}
Our first target has been that of increasing in Scenario 0 the frequency of some train services in Table~\ref{table:current_timetable}, as to get both a larger number of trains a week and a smaller period.
To this aim, we focused on a time span of 48 hours and defined a set $T_{48}$ of 17 trains, where daily trains appear twice, once per day, and trains $t_2, t_3, t_5, t_7$ and $t_8$ appear just once (for reference, see days Tuesday and Wednesday in Table~\ref{table:current_timetable}).
Then we defined the set of candidate trains $T'_{48}$ comprising trains $t_3$ and $t_7$ on the first day, and trains $t_2, t_5, t_8$ on the second day.

Hence we defined an instance of the {\sc pysp} as $I_0 = (S_0, M, T_{48}, T'_{48}, 48 h)$.
Note that in defining this instance, we added Constraints~\eqref{eq:unavailability_constraints} and~\eqref{eq:unavailability_activation_constraints} (see \ref{appendix:unavail}) to the {\sc milp} model, in order to model the unavailability of the reach stackers during the night shift.
In 776 seconds -- all computations were run on a laptop with a 11th Gen Intel Core i7 @ 2.80GHz × 8 processor using 16 GB of RAM and running Debian 11, using Cplex 20.1. -- we found the optimal solution, entailing 18 trains every 48 hours, i.e. 54 trains a week: the solution is such that $t_2, t_3, t_5$ and $t_7$ are operated three times a week while $t_8$ and the other train services are operated daily. The resulting traffic volume corresponds to a 15$\%$ increment with respect to the weekly baseline of 47 trains.

\smallskip
In order to possibly increase the number of weekly trains we moved to Scenario 1. We first solved an instance analogous to the one of the previous scenario, i.e. we defined $I_1' = (S_1, M, T_{48}, T'_{48}, 48 h)$ but with the reach stackers operated 24 hours a day. The optimal solution of this instance, found in 20805 seconds, entails $22$ trains every 48 hours, i.e. every train in $T_{48} \cup T'_{48}$. This proved to us that in Scenario 1 it is possible to operate a timetable with a period of 24 hours.

At this point we defined yet another instance to see whether: $(i)$ the same timetable could be operated with a schedule with period 24 hours, $(ii)$ new train services from a set $T'$ suggested by RFI could be added. The set $T'$ is made of 5 additional sevices for which RFI suggested time windows for the arrival time and for the departure as to agree with our model.
Hence we defined $I_1 = (S_1, M, T, T', 24 h)$. $I_1$ was solved to optimality in 15636 seconds and its optimal solution operates 12 trains per day, i.e. only one train service from $T'$ could be served for a total number of 72 trains a week, corresponding to a 53$\%$ increment with respect to the weekly baseline.

\subsection{Scenarios 2 and 3}
In order to possibly increase the number of weekly trains we moved to Scenario 2 ad 3. We defined instances $I_2 = (S_2, M, T, T', 24 h)$ and $I_3 = (S_3, M, T, T', 24 h)$. Unfortunately, we could not solve these instance to optimality within the established time limit of 16 hours. After 16 hours of computation the best available feasible solutions were a solution with 13 trains out of 16 for {\em Scenario 2}, and one with 12 trains out of 16 for {\em Scenario 3}. These instances are harder to solve because the higher number of trains than can be scheduled determines a significant expansion of the set of potential conflicts, and, in turn, of the generated constraints. Moreover, in Scenario 3 each train has twice the number of plans available in the previous scenario -- because for each train we need to choose whether using reach stackers or the rail mounted crane, for loading and unloading, further contributing to the complexity of instance $I_3$.

In order to complete the analysis of these scenarios, we employed a heuristic approach based our {\sc milp} and consisting of two steps, as to exploit the solution computed for the previous scenario as a starting point. We provide more detail about this heuristic in Appendix \ref{appendix:heuristic}.

By following this procedure, and starting from the optimal solution for Scenario 1, for both scenarios we were able to find solutions operating 16 trains, i.e. all trains in $T \cup T'$, in 337 and 31 seconds respectively, proving that under the hypotheses of Scenarios 2 and 3 it is possible to handle 96 trains a week, corresponding to a 104$\%$ increment with respect to the weekly baseline.
Observe that, since the number of trains scheduled matches the total number of candidate trains, these solutions are optimal even if obtained through a heuristic procedure.

\begin{figure*}
        \centering
        \begin{subfigure}[b]{0.475\textwidth}
            \centering
            \includegraphics[width=\textwidth]{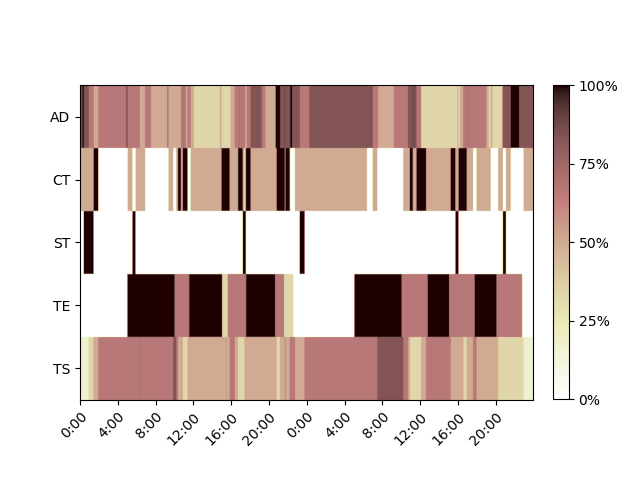}
            \caption{Scenario 0: 18 trains every 48 hours.}
            \label{fig:heatmap_0}
        \end{subfigure}
        \hfill
        \begin{subfigure}[b]{0.475\textwidth}  
            \centering 
            \includegraphics[width=\textwidth]{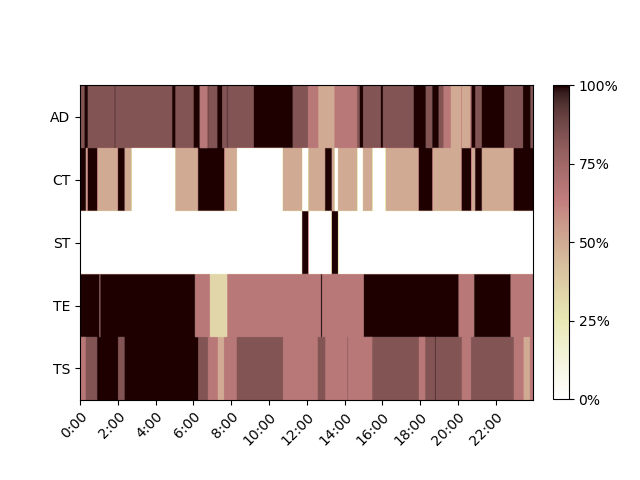}
            \caption{Scenario 1: 12 trains every 24 hours.}
            \label{fig:heatmap_1}
        \end{subfigure}
        \vskip\baselineskip
        \begin{subfigure}[b]{0.475\textwidth}   
            \centering 
            \includegraphics[width=\textwidth]{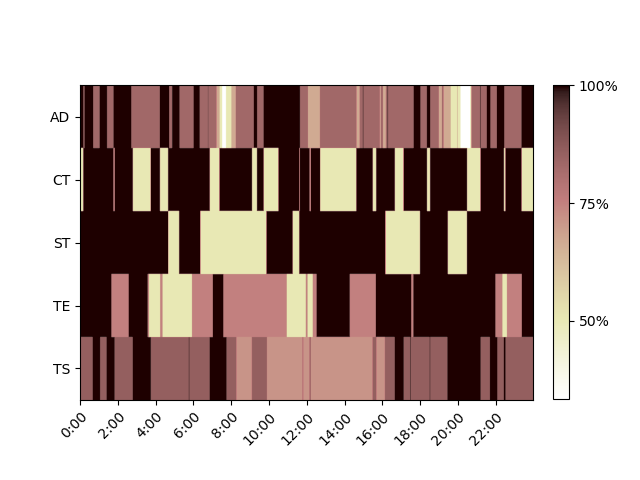}
            \caption{Scenario 2: 16 trains every 24 hours.}
            \label{fig:heatmap_2}
        \end{subfigure}
        \hfill
        \begin{subfigure}[b]{0.475\textwidth}   
            \centering 
            \includegraphics[width=\textwidth]{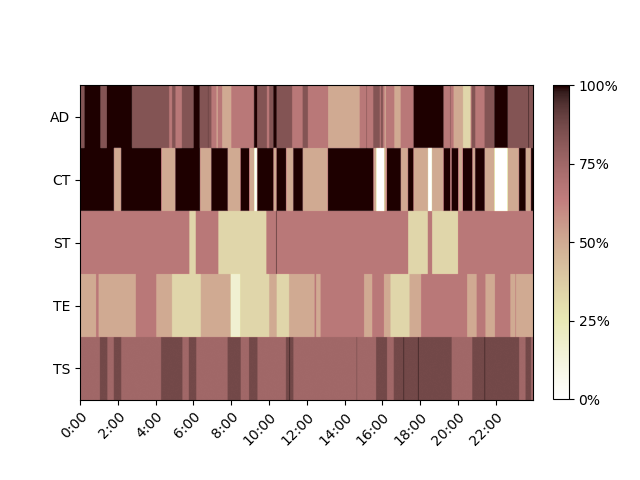}  
            \caption{Scenario 3: 16 trains every 24 hours.}
            \label{fig:heatmap_3}
        \end{subfigure}
        \caption{Heatmaps representing the instantaneous utilization over time for different resources (i.e. arrival/departure tracks, connection tracks, side tracks, transshipment equipment, transshipment equipment), computed on the schedules obtained on each of the four different scenarios.} 
        \label{fig:heatmaps}
    \end{figure*}

\smallskip
In order to better analyze the solutions computed for the different scenarios, we show in Figure~\ref{fig:heatmaps} the resource utilization in each scenario. According to Figure~\ref{fig:heatmap_0}, the current bottleneck for the transshipment yard of Marzaglia is the transshipment equipment, i.e. the reach stackers: for 58$\%$ of their work shift there are no reach stackers available and their average utilization is slightly less than 85$\%$. With the exception of reach stackers, the average utilization of all resources is quite modest and the obtained solution is overall easily and safely implemented in practice.
The reach stackers are the critical resource also in Scenario 1, as can be seen in Figure~\ref{fig:heatmap_1}: their average utilization is $83\%$ and there are no available reach stackers for 54$\%$ of the day. Even if still implementable in practice, this schedule determines a greater congestion with respect to the previous scenario.

We now move to Scenario 2. In this case, identifying the bottleneck is not that easy, as it is in  illustrated in Figure~\ref{fig:heatmap_2}. Here, in addition to the transshipment equipment, the connection tracks and the side tracks are also close to saturation: these resources are completely saturated respectively for $49\%$, $70\%$ and $70\%$ of the day. Likely, increasing the capacity of just one of these resources would not be sufficient to increase the practical capacity of the yard, also because every resource of the yard is indeed quite congested as their average utilization of each resource falls in the range 80-85$\%$. In this context that most areas of the yard are so intensively utilized because they are also used by trains as {\em waiting areas}, notably the transshipment area, the arrival/departure area, and of course the side track area. This is a sign that the transshipment yard suffers from a non-ideal timetable, imposed by the transit line; indeed, the average stay of a train, according to the given timetable, is about twice as long as the time required for traversing the yard and completing the transshipment operations.

Finally Figure~\ref{fig:heatmap_3} shows that in Scenario 3 the arrival/departure area, the connection tracks and the transshipment tracks all have a high utilization (ranging from $76\%$ to $85\%$). However, only the connection tracks are completely saturated for a longer fraction of the day ($59\%$). By comparing this Figure with Figure $\ref{fig:heatmap_2}$, it is clear that even if both scenarios can operate twice the current weekly number of trains, if it is likely to double the current traffic volume, Scenario 3 is a much safer choice, even if more expensive.

\section{Conclusions}
In this paper we proposed an approach for the capacity evaluation of a transshipment yard, that we have developed for Rete Ferroviaria Italiana (RFI), the Italian railway infrastructure manager.
Our approach allows to analyze the practical capacity of transshipment yards, keeping into account all the possible aspects related to layout and timetable constraints. The tool has been developed for analyzing multiple transshipment yards of the RFI network, identifying their current bottlenecks and evaluating the impact of different investment plans finalized at boosting their capacity. In this paper we reported one of this case studies: the transshipment yard of Marzaglia. 
Our approach is based on what we called the Yard Saturation Problem, consisting in saturating a (possibly empty) initial timetable, providing not only an estimate of the practical capacity of a transshipment yards, useful in strategic planning, but also operational instructions for operating the yard at its full capacity.
Indeed, by solving the Periodic version of the Yard Saturation Problem a periodic schedule can be obtained for all operations of each train, safely and easily implementable in practice. 

Through the application of our approach to the transshipment yard of Marzaglia we were able to provide a reliable estimate of the current capacity and bottlenecks, together with an evaluation of three different investment plans; in particular, we showed that each investment scenario would allow to operate each of the current train services on a daily basis, while the two scenarios with higher cost further allow doubling the number of served trains per week.

RFI is determined to apply our method to several other transshipment yards of its network, and while this can be carried out right away for medium sized transshipment yard, we plan to improve the scalability of our model in order to tackle larger yards. A first heuristic resolution scheme has been already tested during the Marzaglia case study, but the design of more complex heuristic algorithms and/or decomposition strategy will be object of future work.

\printbibliography

\appendix

\section{Alternative graph}

\subsection{Periodic Yard Saturation Problem}
\label{appendix:periodic_scheduling}
In the following we will assume that only a single plan is available for each train, i.e. $P(t) = \{p^t\} ~ \forall t \in T \cup T'$.
In order to include the periodic constraints in the alternative graph representation that we have described in Section~\ref{section:disjunctive_graphs}, for each conflicting operations pair $\{u, v\} \in F$, we now have to include in $A_c$ a set $A(\{u, v\})$ made of $4(2k - 1)$ arcs, i.e. $2k-1$ tuples, consisting of two precedence arcs and two meeting arcs, each handling the conflict between $u$ and $v^i$, or $u^i$ and $v$ for $i=1..k$.

Let us consider one of these tuples defined for the pair of replicas $(u, v^i), i = 1, \ldots, k$ (pairs $(v, u^i)$ are handled analogously).
We let $u'$ be the successor of $u$ and $v'$ be the successor of $v$. Then we add to the arc set $A_c$:
\begin{itemize}
\item a precedence arc for $u$ over $v^i$: $(u', v)$, having length $\epsilon + (i-1)\tau$;
\item a precedence arc for $v^i$ over $u$: $(v', u)$, having length $\epsilon - (i-1)\tau$;
\item meeting arcs: $(u, v')$, having length $-\epsilon + (i-1)\tau$, and $(v, u')$, having length $-\epsilon - (i-1)\tau$.
\end{itemize}

\subsection{Multiple plans for each train service}
\label{appendix:multiple_plans}

Here we provide additional details on how the alternative graph should be adjusted when multiple plans are available for each train service.
As explained in Section~\ref{section:multiple_plans}, we consider an alternative graph where, for any train $t$, each arc in $A(t)$ is associated to one or more plans, and the set of vertices $V(t)$ induces a subgraph where $(i)$ arcs in $\overrightarrow{A}(t)$ form an acyclic graph and $(ii)$ each plan in $P(t)$ corresponds to an oriented path in this subgraph, from the arrival node to the departure node (see an example in Figure~\ref{fig:multiple_plans}).

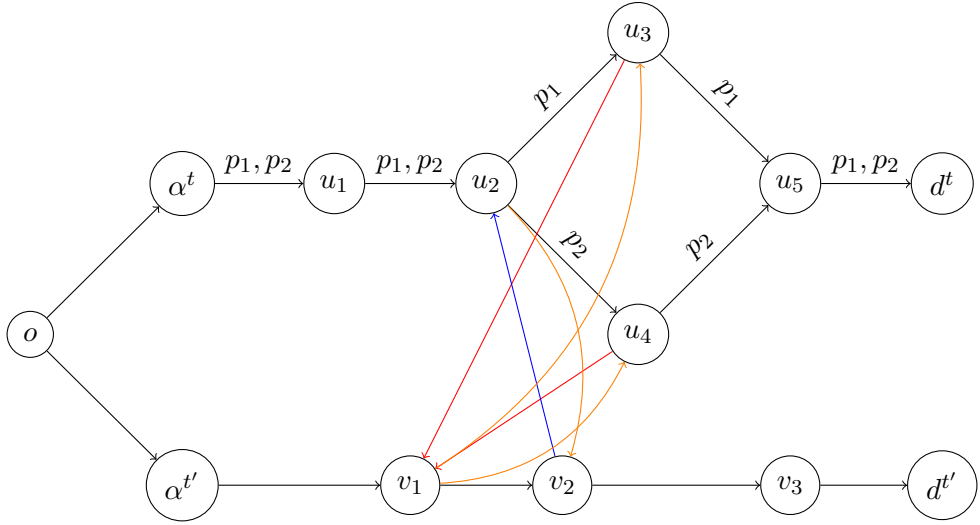
\begin{figure}
\begin{center}
\begin{tikzpicture}[main/.style = {draw, circle}] 
\pgfmathsetmacro{\l}{2}
\node[main] (1) at (0*\l,-2*\l) {$o$};
\node[main] (8) at (1*\l,-1*\l) {$\alpha^t$};
\node[main] (2) at (2*\l,-1*\l) {$u_1$};
\node[main] (3) at (3*\l,-1*\l) {$u_2$};
\node[main] (4) at (4*\l,0*\l) {$u_3$};
\node[main] (5) at (4*\l,-2*\l) {$u_4$};
\node[main] (6) at (5*\l,-1*\l) {$u_5$};
\node[main] (7) at (6*\l, -1*\l) {$d^t$};
\node[main] (9) at (1*\l,-3*\l) {$\alpha^{t'}$};
\node[main] (10) at (2.5*\l,-3*\l) {$v_1$};
\node[main] (11) at (3.5*\l,-3*\l) {$v_2$};
\node[main] (12) at (5*\l,-3*\l) {$v_3$};
\node[main] (13) at (6*\l, -3*\l) {$d^{t'}$};

\begin{scope}[every edge/.style={draw=black}]
\path [->] (1) edge node [pos=0.5, sloped, above]{} (8);
\path [->] (8) edge node [pos=0.5, sloped, above]{$p_1, p_2$} (2);
\path [->] (2) edge node [pos=0.5, sloped, above]{$p_1, p_2$} (3);
\path [->] (3) edge node [pos=0.5, sloped, above]{$p_1$} (4);
\path [->] (3) edge node [pos=0.5, sloped, above]{$p_2$} (5);
\path [->] (4) edge node [pos=0.5, sloped, above]{$p_1$} (6);
\path [->] (5) edge node [pos=0.5, sloped, above]{$p_2$} (6);
\path [->] (6) edge node [pos=0.5, above]{$p_1, p_2$} (7);
\path [->] (1) edge node [pos=0.5, sloped, above]{} (9);
\path [->] (9) edge node [pos=0.5, sloped, above]{} (10);
\path [->] (10) edge node [pos=0.5, sloped, above]{} (11);
\path [->] (11) edge node [pos=0.5, sloped, above]{} (12);
\path [->] (12) edge node [pos=0.5, sloped, above]{} (13);
\end{scope}
\begin{scope}[every edge/.style={draw=red}]
\path [->] (4) edge node [red, pos=0.5, sloped, above]{} (10);
\path [->] (5) edge node [red, pos=0.5, sloped, above]{} (10);
\end{scope}
\begin{scope}[every edge/.style={draw=blue}]
\path [->] (11) edge node [blue, pos=0.5, sloped, above]{} (3);
\end{scope}
\begin{scope}[every edge/.style={draw=orange}]
\path [->] (3) edge [bend left] node [orange, pos=0.5, sloped, above]{} (11);
\path [->] (10) edge [bend right] node [orange, pos=0.5, sloped, above]{} (4);
\path [->] (10) edge [bend right] node [orange, pos=0.5, sloped, above]{} (5);
\end{scope}
\end{tikzpicture}
\end{center}
\caption{Precedence arcs (in red and blue) and meeting arcs (orange) for the conflicting pair of operations $\{u_2, v_1\}$ belonging to trains $t$ and $t'$. Note that for simplicity we are omitting some arcs in $A(t)$ and $A(t')$, as well as the length of the arcs in $A_{c}$.}
\label{fig:multiple_plans_conflict}
\end{figure}

An immediate consequence of this extension is that now each node $u \in V(t), t \in T \cup T'$  may have multiple successors, depending on the chosen plan.
Therefore we introduce the set $\Sigma(u)$ to indicate the set of all successors of node $u$.
More precisely $\Sigma(u) = \{v: \exists ~ p^t \in P(t), t \in T \cup T', i \geq 0 ~ s.t. ~ u = p^t(i), v = p^t(i+1)\}$.

For each conflicting operations pair $\{u, v\} \in F$, as before we define $2k-1$ tuples of arcs in $A(\{u, v\}) \subseteq A_c$. However, each tuple now contains $|\Sigma(u)| + |\Sigma(v)|$ precedence arcs and $|\Sigma(u)| + |\Sigma(v)|$ meeting arcs. More in depth, let us consider the conflict associated with the pair of replicas $(u ,v^i), i = 1, \ldots, k$.
With this conflict we now associate the precedence arcs:
\begin{itemize}
\item $(u', v)$, having length $\epsilon + (i-1)\tau$, $\forall u' \in \Sigma(u)$;
\item $(v' u)$, having length $\epsilon - (i-1)\tau$, $\forall v' \in \Sigma(v)$;
\end{itemize}
and the meeting arcs:
\begin{itemize}
\item $(v, u')$, having length $-\epsilon - (i-1)\tau$, $\forall u' \in \Sigma(u)$;
\item $(u, v')$, having length $-\epsilon + (i-1)\tau$, $\forall v' \in \Sigma(v)$;
\end{itemize}
Figure~\ref{fig:multiple_plans_conflict} shows an example of one of such tuples, where one of the nodes involved in the conflict has more than one successor.

\section{{\sc milp}-formulation}
Here we provide additional technical details on the constraints of the {\sc milp}-model of the {\sc PYSP}, and on the resolution scheme that we adopt.

\subsection{Multiple plans for each train}
\label{appendix:multiple_plans_milp}

As stated in \ref{appendix:multiple_plans}, when multiple plans are available for each train, each conflict between two conflicting replicas $u$ and $v^i, i = 1..k$, is modeled through $|\Sigma(u)|$ precedence arcs for $u$, $|\Sigma(v)|$ precedence arcs for $v^i$, and $|\Sigma(u)| + |\Sigma(v)|$ meeting arcs.
However, in any feasible solution to the {\sc PYSP} we want exactly one precedence arc or two meeting arcs to be active, analogously to what happens for the single plan case.
In particular, supposing that we want to give precedence to $u$ over $v^i$ we want to activate, among the $|\Sigma(u)|$ precedence arcs, the precedence arc that is associated with the active successor of $u$ in $\Sigma(u)$.
For this reason, the activation variables $x_a$ for $a \in A_c$ no longer coincide with the $y$ and $z$-variables defined for each conflict. Instead we have to enforce the following constraints.

Let $\{u, v\} \in F$ be a pair of conflicting operations and let us consider the tuple of arcs associated with nodes $u$ and $v^i, i = 1..k$.
For each arc $a = (u', v), u' \in \Sigma(u)$ we enforce:
\begin{equation}
x_{a} \geq x_{uu'} + y_{uv^i} - 1
\end{equation}
For each arc $a = (v', u), v' \in \Sigma(v)$ we enforce:
\begin{equation}
x_{a} \geq x_{vv'} + y_{v^iu} - 1
\end{equation}
For each arc $a = (u, v'), v' \in \Sigma(v)$ and $a' = (v, u'), u' \in \Sigma(u)$ we enforce, respectively:
\begin{equation}
x_{a} \geq x_{vv'} + z_{\{uv^i\}} - 1
\end{equation}
\begin{equation}
x_{a'} \geq x_{uu'} + z_{\{uv^i\}} - 1
\end{equation}

\subsection{Capacity cuts}
\label{appendix:capacity_cuts}

Here we provide additional details on how the upper bound on $z^*$ of Inequalities~\eqref{eq:cap_risorse} and~\eqref{eq:cap_risorsebis} in Section~\ref{section:analytical_bounds} can be linearized and added to the {\sc MILP}-formulation.
In particular, from Inequality~\eqref{eq:cap_risorse} we derive the following cutting plane:
\begin{equation}
\label{eq:cut_1}
\sum_{t \in T \cup T'} \underline{\rho}_t \phi_t \leq b \tau
\end{equation}
while, starting from~\eqref{eq:cap_risorsebis}, we can define the following cut:
\begin{equation}
\label{eq:cut_2}
\sum_{t \in T \cup T'} \phi_t \leq min_{s \in S} \lfloor \frac{u_s \tau}{min_{t \in T \cup T'} \mu_s^t} \rfloor
\end{equation}
Recall that $\mu^t_s$ is the minimum occupation time of $s$ by $t$.
Moreover, from Inequality~\eqref{eq:cap_risorsebis}, we can also derive the following family of cuts, where $\mu^{t, p^t}_s$ is the minimum occupation time of $s$ by $t$ if it follows plan $p^t$:
\begin{equation}
\label{eq:cut_3}
\sum_{t \in T \cup T'} \sum_{p^t \in P(t)} w_{p^t} \mu_s^{t, p^t} \leq \tau u_s \qquad s \in S
\end{equation}

\subsection{Bound on average utilization}
\label{appendix:robustness}

In order to make sure that the average utilization of any resource is not higher than $85\%$ of the full capacity, we add the following sets of variables and constraints.
First we define a continuous variable $g^s_{p^t} \in \mathbb{R}_{\geq 0}$ for each resource $s \in S$, train $t \in T \cup T'$ and plan $p^t \in P(t)$, representing the total occupation time of resource $s$ by train $t$ if $p^t$ is active.
Then we add the following constraints to activate these variables:
\begin{equation}
\label{eq:constraint_utilization_1}
g^s_{p^t} \geq \sum_{j \in [n_{p^t} - 1] : s \in S_{p^t(j)}} (\sigma(p^t(j+1)) - \sigma(p^t(j))) + \epsilon \beta^s_{p^t} - M (1 - w_{p^t})
\end{equation}
where $\beta^s_{p^t}$ is the number of times that train $t$ acquires and releases resource $s$.
Finally, we add a constraint for each resource, ensuring that the total average utilization throughout a period does not exceed the given threshold of $85\%$:
\begin{equation}
\label{eq:constraint_utilization_2}
\sum_{t \in T \cup T'} \sum_{p^t \in P(t)} g^s_{t, p^t} \leq 0.85 \cdot u_s \tau \qquad s \in S
\end{equation}

\smallskip

Note that this constraint allows to strengthen the capacity cuts defined by Inequalities~\eqref{eq:cut_1}-\eqref{eq:cut_3}, by multiplying the right hand side by the constant $0.85$. This is shown in the complete {\sc milp}-formulation in \ref{appendix:complete_milp}.

\subsection{Resource unavailabilities}
\label{appendix:unavail}
Here we describe how the alternative graph $G(V, A, l)$ of Section~\ref{section:disjunctive_graphs}, and the {\sc milp}-formulation of Section~\ref{section:rfi_milp} should be adjusted to take into account temporary resource unavailabilities.

Let us consider a resource $s$ that is unavailable from $h$ to $h'$, with $[h, h') \subset [0, \tau)$.
A simple way to model this resource unavailability, used in \cite{burdett_2009}, is by introducing $u_s$ artificial trains in $T$, which must use resource $s$ in $[h, h')$.
One drawback of that method is that it forces every other real operation requiring $s$ to be either completed before $h$ or started after $h'$.
This can be a very limiting modeling choice for time-expensive operations, as for the transshipment operations, which in practice 
can be started before the interruption forced by the work shifts and completed at the end of the unavailability time interval.
Because of this reason, instead of introducing fictitious trains, we explicitly represent through some additional arcs the three possible choices: a) the operation is completed before $h$, b) started after $h'$, c) started before $h$ and completed after $h'$.

\begin{figure}
\begin{center}
    \begin{tikzpicture}[main/.style = {draw, circle, minimum size=1cm}] 
\pgfmathsetmacro{\l}{2}
\node[main] (1) at (0*\l,-1*\l) {$o$};
\node[main] (2) at (1*\l,-1*\l) {$\alpha^t$};
\node[main] (3) at (2*\l,-1*\l) {$u$};
\node[main] (4) at (3*\l,-1*\l) {$v$};
\node[main] (5) at (4*\l,-1*\l) {$d^t$};

\begin{scope}[every edge/.style={draw=blue}]
\path [->] (1) edge [bend left] node [blue, pos=0.5, sloped, above]{$h'$} (3);
\end{scope}
\begin{scope}[every edge/.style={draw=green}]
\path [->] (4) edge [bend left] node [green, pos=0.5, sloped, above]{$-h$} (1);
\end{scope}
\begin{scope}[every edge/.style={draw=red}]
\path [->] (3) edge [bend left] node [red, pos=0.5, sloped, above]{$\lambda_u + h' - h$} (4);
\end{scope}
\begin{scope}[every edge/.style={draw=black}]
\path [->] (1) edge node [pos=0.5, sloped, above] {} (2);
\path [->] (2) edge node [pos=0.5, sloped, above] {} (3);
\path [->] (3) edge node [pos=0.5, sloped, below] {$\lambda_u$} (4);
\path [->] (4) edge node [pos=0.5, sloped, above] {} (5);

\end{scope}

\end{tikzpicture}
\caption{Example representing how to model the unavailability in $[h, h')$ of a resource $s$ required to the operation associated with node $u$.}
\label{fig:unavailability}
\end{center}
\end{figure}
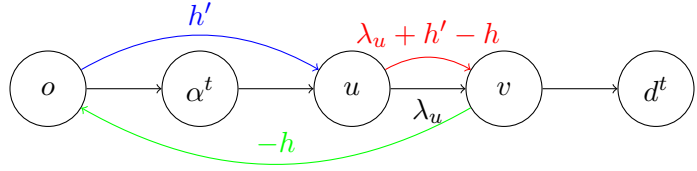

In the alternative graph we add a new class of arcs, that we will denote by $A_{avail}$, including: $(i)$ for each node $v \in \Sigma(u)$, an arc $(v, o)$ of length $-h$, representing case a); $(ii)$ an arc $(o, u)$ of length $h'$, representing case b); $(iii)$ for each node $v \in \Sigma(u)$, an arc $(u, v)$ of length $\lambda_u + (h' - h)$, representing case c).
Figure~\ref{fig:unavailability} shows an example of the arcs modeling the unavailability of a resource.

If we have to consider multiple periods and replicas, i.e. if $k > 1$, then we have to add a larger set of arcs to the graph. Given a node $u$, associated with an operation requiring $s$, for each time window $i = 1, \ldots, k$ the set $A_{avail}$ must include:
\begin{itemize}
\item for each node $v \in \Sigma(u)$, an arc $(v, o)$ of length $-(h + (i-1)\tau)$, representing the case in which the operation associated with $u$ is completed before $h + (i-1)\tau$; we will refer to this arc as $a^{i, before}_{u, v}$;
\item arc $(o, u)$ of length $h' + (i-1)\tau$, representing the case in which the operation associated with $u$ starts after $h' + (i-1)\tau$; we will refer to this arc as $a^{i, after}_{u}$;
\item for each node $v \in \Sigma(u)$, an arc $(u, v)$ of length $\lambda_u + (h' - h)$, representing the case in which the operation associated with $u$ is temporarily interrupted by the unavailability of some resource; we will refer to this arc as $a^{during}_{u, v}$.
\end{itemize}

\medskip

Regarding the {\sc milp}-formulation, for each operation $u$ requiring $s$, i.e. $s \subseteq S_u$, and for each period $i = 1, \ldots, k$ we introduce three binary variables:
\begin{itemize}
	\item $\underline{r}_u^i$, that is set to 1 if operation $u$ is completed before the start of the time interval in period $i$;
	\item $\overline{r}_u^i$, that is set to 1 if operation $u$ is started after the end of the time interval in period $i$;
	\item $r_u^i$, that is set to 1 if operation $u$ is interrupted by the time interval.
\end{itemize}
Since we have to choose exactly one of these three alternatives for each active node and each period, we add the following \textit{unavailability constraints}:
\begin{equation}
\label{eq:unavailability_constraints}
\underline{r}_u^i + \overline{r}_u^i + r_u^i = x_u \qquad s \in S_u, i = 1, ..., k
\end{equation}

Moreover, we have to associate an arc activation variable $x_a$ to each arc $a \in A_{avail}$, as we have done with all the other arcs in $A$. In particular, for each operation $u$ requiring $s$ we enforce the following arc activation constraints:
\begin{align}
\label{eq:unavailability_activation_constraints}
x_{a^{i, before}_{u, v}} & \geq x_{uv} + \underline{r}^i_u - 1 & v \in \Sigma(u), i = 1, \ldots k\\
x_{a^{i, after}_{u}} & \geq \overline{r}^i_u & i = 1, \ldots k\\
x_{a^{during}_{u, v}} & \geq x_{uv} + \sum_{i=1}^{k} r_u^i - 1 & v \in \Sigma(u)
\end{align}

\bigskip

In a more general setting, we can allow multiple intervals of unavailabilities for $s$, i.e. $[h_1, h_1'), \ldots, [h_d, h_d') \subset [0, \tau)$.
However, we will assume that the maximum waiting time $\gamma_u$ of each operation $u$ requiring $s$ does not allow more than one interruption per period. This assumption is not strictly necessary but is compatible with our case study and greatly reduces the technical complexity of the graph representation. Indeed, under this assumption, we can manage multiple intervals of unavailabilities independently from each other.

\subsection{Complete {\sc milp}-model}
\label{appendix:complete_milp}

The complete {\sc milp}-model is:
\begin{flalign*}
&max \sum_{t \in T \cup T'} \phi_t &
\end{flalign*}
\begin{flalign*}
&\phi_t = 1 & t \in T \\
&\sum_{p \in P(t)} w_{p} = \phi_t & t \in T \cup T'
\end{flalign*}
\begin{flalign*}
&x_u = \begin{cases}
        1 \qquad \qquad \quad ~ \text{if $u = o$,}\\
        \sum_{p \in P(t): u \in V(p)} w_{p} ~~ \text{otherwise}
        \end{cases} & u \in V
\end{flalign*}
\begin{flalign*}
&x_a = \begin{cases}
        \phi_t  ~~ \text{if $a \in A(t)_{a} \cup A(t)_{d}$,}\\
        \sum_{p \in P(t): a \in A(p)} w_{p} ~~ \text{otherwise}
        \end{cases} & a \in A(t), t \in T \cup T'
\end{flalign*}
\begin{flalign*}
\begin{split}
    y_{uv^i} + y_{v^iu} + z_{\{uv^i\}} \geq x_u + x_v - 1 \\
    y_{vu^i} + y_{u^iv} + z_{\{vu^i\}} \geq x_u + x_v - 1
\end{split}
& \{u, v\}\in F, i= 1..k
\end{flalign*}
\begin{flalign*}
\begin{split}
x_{a} \geq x_{uu'} + y_{uv^i} - 1 \\
x_{a} \geq x_{uu'} + y_{u^iv} - 1 \\
\end{split} & \{u, v\} \in F, a \in A(\{u, v\}), i=1..k, u' \in \Sigma(u) \text{ tail of } a, v \text{ head of } a
\end{flalign*}
\begin{flalign*}
\begin{split}
x_{a} \geq x_{vv'} + y_{v^iu} - 1 \\
x_{a} \geq x_{vv'} + y_{vu^i} - 1 \\
\end{split} & \{u, v\} \in F, a \in A(\{u, v\}), i=1..k, v' \in \Sigma(v) \text{ tail of } a, u \text{ head of } a
\end{flalign*}
\begin{flalign*}
\begin{split}
  x_{a} \geq x_{vv'} + z_{\{uv^i\}} - 1 \\
  x_{a} \geq x_{vv'} + z_{\{u^iv\}} - 1 \\  
\end{split} & \{u, v\} \in F, a \in A(\{u, v\}), i=1..k, v' \in \Sigma(v) \text{ head of } a, u \text{ tail of } a
\end{flalign*}
\begin{flalign*}
\begin{split}
    x_{a} \geq x_{uu'} + z_{\{uv^i\}} - 1 \\
    x_{a} \geq x_{uu'} + z_{\{u^iv\}} - 1 \\
\end{split} & \{u, v\} \in F, a \in A(\{u, v\}), i=1..k, u' \in \Sigma(u) \text{ head of } a, v \text{ tail of } a
\end{flalign*}
\begin{flalign*}
& \sum_{\{u, v^i\} \subseteq Q} z_{\{uv^i\}} \leq \binom{u_s + 1}{2} - 1 & \forall s \in S, \forall Q \subseteq O_s, |Q| = u_s + 1
\end{flalign*}
\begin{flalign*}
& \sigma(v) - \sigma(u) \geq l_{a} - M(1 - x_{a}) & a \in A, v \text{ tail of } a, u \text{ head of } a \\
& \sigma(o) = 0 &
\end{flalign*}
\begin{flalign*}
&\sum_{t \in T \cup T'} \phi_t \leq min_{s \in S} \lfloor \frac{0.85 \cdot u_s \tau}{min_{t \in T \cup T'} \mu_s^t} \rfloor &
\end{flalign*}
\begin{flalign*}
& \sum_{t \in T \cup T'} \sum_{p \in P(t)} w_{p} \mu_s^{t, p} \leq 0.85 \cdot \tau u_s & s \in S
\end{flalign*}
\begin{flalign*}
& \sum_{t \in T \cup T'} \underline{\rho}_t \phi_t \leq 0.85 \cdot b \tau &
\end{flalign*}
\begin{flalign*}
& g^s_{p^t} \geq \sum_{j \in [n_{p^t} - 1] : s \in S_{p^t(j)}} (\sigma(p^t(j+1)) - \sigma(p^t(j))) + \epsilon \beta^s_{p^t} - M (1 - w_{p^t}) & s \in S, t \in T \cup T', p^t \in P(t) \\
& \sum_{t \in T \cup T'} \sum_{p^t \in P(t)} g^s_{t, p^t} \leq 0.85 \cdot u_s \tau & s \in S
\end{flalign*}
$$ \phi \in \{0, 1\}^{|T|}, \qquad w \in \{0, 1\}^{|P|}, \qquad x \in \{0, 1\}^{|V| + |A|} $$
$$ y \in \{0, 1\}^{2|F|(2k -1)}, \qquad z \in \{0, 1\}^{|F|(2k-1)}, \qquad \sigma \in \mathbb{R}_{\geq 0}^{|V|}, \qquad g \in \mathbb{R}_{\geq 0}^{|S||P|}$$

\bigskip

Note that in presence of a temporary unavailability of some resource, we may need to add other suitable variables and constraints that are described in \ref{appendix:unavail}.

\bigskip

Observe that for arc $a = (u, v) \in A$ a safe choice for constant $M$ is $M = ub(u) - lb(v) + l_{uv}$, where $ub(u)$ is a known upper bound for $\sigma(u)$ and $lb(v)$ is a known lower bound on $\sigma(v)$. For the origin node we can take $lb(o) = ub(o) = 0$, and for a node $u$ in $V(t)$ an obvious choice is $lb(u) = \underline{a}_t$ and $ub(u) = \overline{q}_t$. However, it is possible to compute tighter bounds by exploiting the information on the minimum and maximum duration $\lambda_m$ and $\gamma_m$ of each operation $m$ involved in the plans available to train $t$.

\subsection{Solving the {\sc MILP}-model}
\label{appendix:algorithm}

In this section we provide additional details on the dynamic row and column generation scheme adopted to solve the {\sc MILP}-model of the {\sc PYSP}.

We define an iterative scheme that in the generic iteration $h$ works with:
\begin{itemize}
	\item a restricted set of conflicting operations pairs $F_h$;
	\item a restricted alternative graph $G_h(V, A_h, l_h)$, where $A_h = A_s \cup A_c^h$;
	\item a set $\mathcal{Q}_h^s$ for each resource $s$, containing sets of size $u_s + 1$ of operations requiring $s$.
\end{itemize}
In each iteration a {\sc MILP} problem $R_h$, based on the alternative graph $G_h$ and the set $F_h$, $\mathcal{Q}_h^s, s \in S$, is solved.
The obtained solution provides set of schedules $\mathcal{P}^h = \{(p^t_h, \sigma^t_h), t \in \bar{T}\}$ for some $T \subseteq \bar{T} \subseteq T \cup T'$.
This set of schedules is, in general, not feasible, as it is obtained by enforcing just a restricted set of constraints.

In particular, at the start of the algorithm we have:
\begin{itemize}
	\item $F_0 = \emptyset$;
	\item $A_c^0 = \emptyset$;
	\item $\mathcal{Q}_0^s = \emptyset$ for each resource $s$.
\end{itemize}

\paragraph{Row generation}
At each iteration the set of schedules $\mathcal{P}^h$ is checked for any capacity violation on each resource in $S$.
This can be done by constructing an interval graph for each resource $s$ and by checking it for cliques $Q^s$ of size $u_s + 1$.
The nodes of the interval graph for resource $s$ are the active nodes $u$ in $G_h$ that are associated with operations requiring $s$, and all their replicas $u^i$, $i = 2, ..., k$.
Let us consider an active node $u$, we denote by $u'$ the unique active node in $\Sigma(u)$.
The time interval $\mu(u)$ associated with node $u$ (and its operation), according to the schedule $\mathcal{P}^h$, is $\mu(u) = [\sigma(u), \sigma(u') + \epsilon)$.
For each replica $u^i, i = 2, ..., k$, the associated time interval is $\mu(u^i) = [\sigma(u) + (i-1)\tau, \sigma(u') + (i-1)\tau + \epsilon)$.
In the interval graph two nodes $u^i$ and $v^j$ are adjacent to each other if $\mu(u^i) \cap \mu(v^j) \neq \emptyset$.
Let $K_h^s$ the set of all maximal cliques on the interval graph for resource $s$.
It is well known that such set can be found in linear time on interval graphs \cite{rolf_1986}.
Then we define:
\begin{itemize}
	\item a new set of conflicting pairs of operations $F_{h+1}$, by adding all the pairs of operations contained in each of the cliques of size greater than $u_s$;
	\item a new set $A_c^{h+1}$ containing all precedence and meeting arcs associated to the pairs in $F_{h+1}$;
	\item new sets $\mathcal{Q}_{h+1}^s$ which now includes the cliques found for each resource $s$, if any.
\end{itemize}
Note that the resulting problem $R_{h+1}$ will differ from $R_h$ both for the constraints set and the variables sets.

If at any point the set $K_h$ is empty, then $\mathcal{P}^h$ is a feasible schedule and the algorithm stops.
Algorithm~\ref{alg:bigM} illustrates the pseudocode for the algorithm.

\begin{algorithm}
   \caption{Big-M formulation algorithm}
    \label{alg:bigM}
    \begin{algorithmic}[1]
      \Function{solve}{$\mathcal{W}, A_0$}
      	\State set of conflicts constraints $\mathcal{K} = \emptyset$
        \State set of capacity constraints $\mathcal{U} = \emptyset$
        \State set of precedence constraints $\mathcal{D}$ defined over $A_0$
        \State set of plan selection constraints $\mathcal{W}$
        \While {true}
            \State $w, x, y, z, \sigma = \textsc{optimize}(\mathcal{K}, \mathcal{U}, \mathcal{P}, \mathcal{W})$
      		\State find set of capacity violations $B$ from $\sigma$
            \If{ $B = \emptyset$}
            		\State \textbf{break}
            	\Else
            		\State update $\mathcal{K}, \mathcal{U}, \mathcal{D}$ based on $B$
             \EndIf
    		\EndWhile
    \EndFunction
\end{algorithmic}
\end{algorithm}

\section{Heuristic procedure}
\label{appendix:heuristic}

Here we describe the two-steps heuristic procedure that we employed for the analysis of Scenarios 2 and 3 in the transshipment yard of Marzaglia (see Section \ref{section:rfi_results}).

Given a schedule ${\cal P} = \{(p^t, \sigma^t), t \in \overline{T}\}$ for a set of trains $\overline{T} \subseteq T \cup T'$, we say that we {\em fix $\cal P$ with tolerance $\eta$}, if we solve an instance of the Saturation Problem where each operation $u$ of a train $t \in \overline{T}$ is bounded to a time interval of size $2\eta$, centered in $\sigma^t$.
Bear in mind that this tolerance $\eta$ refers only to the schedule ${\cal P}$ and does not change the time intervals imposed on the arrival and departure times, i.e. $[\underline{a}_t, \overline{a}_t]$ and $[\underline{q}_t, \overline{q}_t]$ respectively.

We use this kind of fixing in both steps of our heuristic procedure. The goal of the first step is refining the solution computed for the previous scenario so that we can more easily add new trains in the second step. To this aim, said $T^*_{i-1} \subseteq T'$ the set of candidate trains added in the timetable of scenario $i-1$, in scenario $i$ we solve an instance $I^{balance}_i = (S_{i-1}, M, T \cup T^*_{i-1}, \emptyset, \epsilon, 24 h)$, in which we: $(i)$ fix the schedule obtained by solving $I_{i-1}$ with tolerance $120'$; $(ii)$ minimize the maximum average utilization among the resources, considering the augmented capacities of the resource set $S_i$.
In the second step we solve an instance $I^{heur}_i = (S_i, M, T \cup T^*_i, T \cup T', \epsilon, 24 h)$ in which we fix the schedules obtained by solving $I^{balance}_i$ with tolerance $30'$.

\end{document}